\documentclass[a4paper,12pt]{article}
\usepackage[english]{babel}
\usepackage{tikz}
\usetikzlibrary{matrix,arrows,decorations.markings}
\usepackage{amsmath,amsfonts,amssymb,amsthm,url,textcomp}
\usepackage[utf8]{inputenc}
\usepackage{csquotes}
\usepackage{a4wide}
\usepackage[numbers]{natbib}
\usepackage{fancyhdr}
\usepackage{anyfontsize}
\usepackage{hyperref}
\usepackage{hhline}
\usepackage{leftidx}
\usepackage{array}
\usepackage{longtable}
\usepackage{bm}
\usepackage{enumitem}

\allowdisplaybreaks

\numberwithin{theorem}{section}

\numberwithin{theoremm}{subsection}

\numberwithin{theoremmm}{subsubsection}

\theoremstyle{definition}

\newcommand{\Aut}{\operatorname{Aut}}
\newcommand{\Alt}{\operatorname{Alt}}
\newcommand{\PSL}{\operatorname{PSL}}

\newcommand{\Sym}{\operatorname{Sym}}

\newcommand{\Inn}{\operatorname{Inn}}
\newcommand{\id}{\operatorname{id}}

\newcommand{\fix}{\operatorname{fix}}
\newcommand{\Out}{\operatorname{Out}}
\newcommand{\e}{\mathrm{e}}

\newcommand{\M}{\operatorname{M}}

\newcommand{\GL}{\operatorname{GL}}

\newcommand{\PSU}{\operatorname{PSU}}
\newcommand{\PSp}{\operatorname{PSp}}

\newcommand{\Mod}[1]{\ (\textup{mod}\ #1)}

\renewcommand{\k}{\operatorname{k}}

\newcommand{\IN}{\mathbb{N}}

\newcommand{\Sp}{\operatorname{Sp}}

\newcommand{\IF}{\mathbb{F}}

\newcommand{\GU}{\operatorname{GU}}

\newcommand{\Inndiag}{\operatorname{Inndiag}}

\newcommand{\POmega}{\operatorname{P}\Omega}
\newcommand{\GO}{\operatorname{GO}}

\newcommand{\Outdiag}{\operatorname{Outdiag}}

\newcommand{\omicron}{\operatorname{o}}
\newcommand{\q}{\mathfrak{q}}
\newcommand{\Ord}{\operatorname{Ord}}

\newcommand{\ssrm}{\mathrm{ss}}

\newcommand{\Exp}{\operatorname{Exp}}

\newcommand{\sign}{\operatorname{sign}}
\newcommand{\SO}{\operatorname{SO}}

\begin{document}

\title{Documentation for the GAP code file {\tt OrbOrd.txt}}

\author{Alexander Bors, Michael Giudici and Cheryl E. Praeger\thanks{First author's address: Johann Radon Institute for Computational and Applied Mathematics (RICAM), Altenbergerstra{\ss}e 69, 4040 Linz, Austria. E-mail: \href{mailto:alexander.bors@ricam.oeaw.ac.at}{alexander.bors@ricam.oeaw.ac.at} \newline Second and third author's address: The University of Western Australia, Centre for the Mathematics of Symmetry and Computation, 35 Stirling Highway, Crawley 6009, WA, Australia. E-mail: \href{mailto:michael.giudici@uwa.edu.au}{michael.giudici@uwa.edu.au} and \href{mailto:cheryl.praeger@uwa.edu.au}{cheryl.praeger@uwa.edu.au} \newline The first author is supported by the Austrian Science Fund (FWF), project J4072-N32 \enquote{Affine maps on finite groups}. The second and third authors were supported by the Australian Research Council Discovery Project DP160102323. \newline 2010 \emph{Mathematics Subject Classification}: Primary: 20D60. Secondary: 20D05, 20D45. \newline \emph{Key words and phrases:} Finite groups, Automorphism orbits, Element orders, Simple groups, Monster simple group}}

\date{\today}

\maketitle

\abstract{We give a comprehensive description of the functions and variables defined in the authors' GAP code file {\tt OrbOrd.txt}, which serve mainly to compute (bounds on) the number of $\Aut(S)$-orbits on $S$, or the set or number of element orders in $S$ for nonabelian finite simple groups of Lie type $S$.}

\section{Introduction}\label{sec1}

\subsection{What this code documentation is about}\label{subsec1P1}

The GAP code file {\tt OrbOrd.txt} was written during the authors' work on their paper \cite{BGP19a}, and it is available in text form on the first author's homepage under \url{https://alexanderbors.wordpress.com/sourcecode/orbord/source/}. It consists of various defined functions that are useful for computing (bounds on) the number of $\Aut(S)$-orbits on $S$, or the set or number of element orders in $S$, where $S$ is mostly assumed to be a finite simple group of Lie type. The file does, however, also contain an efficient algorithm for computing the number of element orders in $\Sym(n)$, discussed in Section \ref{sec2}. The functions relating to Lie type groups are discussed in Section \ref{sec3}. For each defined function, we describe its syntax and what it outputs, and in many cases, we also give comments pertaining to the idea behind the algorithm or a source in the literature justifying it. There are also a few defined variables, for which we simply describe what their value is.

\subsection{Some notation used throughout this documentation}\label{subsec1P2}

\begin{itemize}
\item $\omega(G)$: the number of $\Aut(G)$-orbits on the finite group $G$;
\item $\Ord(G)$: the set of element orders in the finite group $G$;
\item $\omicron(G)$: the number of distinct element orders in the finite group $G$ (i.e., $\omicron(G)=|\Ord(G)|$);
\item $\epsilon_{\omega}(S):=\frac{\log\log{\omega(S)}}{\log\log{|S|}}$, where $S$ is a nonabelian finite simple group;
\item $\epsilon_{\q}(S):=\frac{\log\log{(\omega(S)/\omicron(S)+3)}}{\log\log{|S|}}$, where $S$ is a nonabelian finite simple group;
\item $O^{p'}(G)$: the subgroup of the finite group $G$ generated by the elements whose order is a power of the prime $p$;
\item Let $G=X_d(\overline{\IF_p})$ be a simple linear algebraic group of adjoint type over the algebraic closure $\overline{\IF_p}$ of the finite (size $p$) field $\IF_p$, and let $\sigma$ be a Lang-Steinberg endomorphism of $G$ (i.e., a surjective algebraic group endomorphism of $G$ with finite fixed point subgroup $G_{\sigma}$). Moreover, let $B$ be any $\sigma$-invariant Borel subgroup of $X_d(\overline{\IF_p})$, and let $T$ be any $\sigma$-invariant maximal torus of $X_d(\overline{\IF_p})$ contained in $B$. We now define two important parameters $t=t(\sigma)$ and $f=f(\sigma)$: $t$ is defined as the unique smallest positive integer (independent of the choice of $B$ and $T$) such that the $t$-th power of the map $\sigma^{\ast}$ on the character group $X(T)$ induced by $\sigma$ is a positive integral multiple of $\id_{X(T)}$, and $f\in\IN^+/2=\{\frac{1}{2},1,\frac{3}{2},\ldots\}$ is such that $\sigma^{\ast}=p^f\sigma_0$ with $\sigma_0^t=\id_{X(T)}$; $f$ also does not depend on the choice of $B$ and $T$. So $p^f=q(\sigma)$ in the notation of \cite{Har92a}, which is also a notation we will be using, and $f\in\IN^+$ unless $X_d\in\{B_2,G_2,F_4\}$ and $t=2$ (in which case $f$ is half of an odd positive integer, and accordingly, $q=\sqrt{p^{2k+1}}$ for some $k\in\IN$). With these parameters $t$ and $f$ now defined, we denote the group $O^{p'}(G_{\sigma})=O^{p'}(\fix(\sigma))$, which is usually a nonabelian finite simple group, by $\leftidx{^t}X_d(q^t)=\leftidx{^t}X_d(p^{ft})$. For example, $A_3(2)=\PSL_4(2)$, $\leftidx{^2}D_6(9)=\POmega^-_{12}(3)$, and $\leftidx{^2}B_2(8)$ is the Suzuki group of order $29120$. This notation is consistently reflected in the syntax of our functions for Lie type groups.
\end{itemize}

\subsection{Tips on how to read this documentation}\label{subsec1P3}

While the documentation itself is rather long, there is lots of repetition, due to the fact that the functions pertaining to groups of Lie type need to be programmed anew for each Lie family. In order for readers to quickly get a good impression of what this code can do for them, it is sufficient to read the following:
\begin{itemize}
\item Section \ref{sec2}, as well as Subsections \ref{subsec3P1} and \ref{subsec3P18}, each of which contains unique functions without analogues in other parts of the documentation. Moreover, note that the functions described in Subsection \ref{subsec2P6} and Subsubsections \ref{subsubsec3P18P4} to \ref{subsubsec3P18P17} were written for very specific checks in proofs found in \cite{BGP19a} and thus probably will not be needed and can be skipped by readers other than the authors of this documentation.
\item Subsections \ref{subsec3P2} (to get an impression of the functions and their syntax for one family of classical Lie type groups) and \ref{subsec3P8} (the same for one family of exceptional Lie type groups, where mostly the syntax is different).
\end{itemize}
When using {\tt OrbOrd.txt} to study finite simple Lie type groups in general, there are some subtle differences between Lie families to take note of, most notably the following:
\begin{itemize}
\item There are no functions of the form {\tt NrConjugacyClassesX} where {\tt X} is one of the strings {\tt "C"}, {\tt "D"}, {\tt "2D"}, {\tt "E6"}, {\tt "2E6"} or {\tt "E7"}; in those cases, one must use {\tt NrConjugacyClassesBoundX}, which only outputs a lower bound on (not the precise value of) the respective conjugacy class number.
\item For classical groups, some functions have an argument called \enquote{quality level}. This is to be chosen from a certain finite set of positive integers and determines what kind of computations for reaching the desired output (a bound on a certain parameter) are chosen. Generally speaking, the higher the quality level, the better the result (i.e., the sharper the bound that is output), but also the more costly the computations are. Functions for exceptional Lie type groups do not have such a quality level, and there are also differences with regard to the available quality levels among classical groups. For example, the functions of the form {\tt NrAutOrbitsBoundX} always have quality level $2$, but sometimes $1$ as well, available, and the functions of the form {\tt NrElementOrdersBoundX} always have quality levels $1$ and $2$, but sometimes $3$ as well, available.
\item The function {\tt NrElementOrdersSharplyDivisibleByPToEBoundX} is only defined when {\tt X} is either {\tt "A"}, {\tt "2A"}, {\tt "D"} or {\tt "2D"}. This is not because it has no analogues when {\tt X} is {\tt "B"} or {\tt "C"}, but because the authors only needed it in those cases for their work on \cite{BGP19a}.
\item The groups $S$ of type $E_8$ are the only exceptional Lie type groups for which there is no result in the literature (known to the authors) that describes the set of element orders of $S$. Hence there are no functions {\tt ElementOrdersE8} and {\tt NrElementOrdersE8}; instead, we have a function called {\tt NrElementOrdersBoundE8} which, when applied to a prime power $q$, outputs an upper bound on the number of element orders in $E_8(q)$. 
\end{itemize}

\section{Functions for computing the number of element orders in \texorpdfstring{$\Sym(n)$}{Sym(n)}}\label{sec2}

Before we start to list the functions (and the one variable) relevant for this section, let us introduce some (nonstandard) notation and terminology. For $k\in\IN^+$, we denote by $p_k$ the $k$-th prime (so for example, $p_1=2$, $p_2=3$ and so on). Moreover, for a prime $p$ and a positive integer $n$, we denote by
\begin{itemize}
\item $r_p(n)$ the number of unordered integer partitions of $n$ into pairwise coprime prime powers greater than $1$ such that the smallest prime base which occurs is $p$;
\item $r(n)$ the number of unordered integer partitions of $n$ into pairwise coprime prime powers greater than $1$ (so that $r(n)=\sum_{p\leqslant n}{r_p(n)}$ where the index $p$ ranges over primes).
\end{itemize}
For each positive integer $n$, we define the \emph{$n$-th partition number matrix} to be the matrix with $n$ rows and $\max\{1,\pi(n)\}$ columns, where
\[
\pi(n)=|\{p\mid p\text{ is prime and }p\leqslant n\}|,
\]
whose $(i,j)$-th entry, for $i\in\{1,\ldots,n\}$ and $j\in\{1,\ldots,\pi(n)\}$, is $r_{p_j}(i)$. So, for example, the first partition number matrix is the $(1\times 1)$-matrix $(0)$, the second is the $(2\times 1)$-matrix ${0 \choose 1}$, and the third is the $(3\times 2)$-matrix
\[
\begin{pmatrix}0 & 0 \\ 1 & 0 \\ 0 & 1\end{pmatrix}.
\]

\subsection{{\tt PrimesExtended}}\label{subsec2P1}

This is a variable, defined as a comprehensive list of the primes up to $25013$ (the smallest prime larger than $25000$).

\subsection{{\tt NextPartitionNumberMatrix}}\label{subsec2P2}

\begin{itemize}
\item Syntax: {\tt NextPartitionNumberMatrix(m)}, $m$ being the $n$-th partition number matrix (in the above defined sense) for some $n\in\IN^+$.
\item Output: the $(n+1)$-th partition number matrix.
\item Comments: This algorithm uses the fact that the numbers $r_p(k)$ satisfy the recursion
\[
r_p(k)=\sum_{e=1}^{\lfloor\log_p(k)\rfloor}{\begin{cases}\sum_{p<\ell\leqslant k,\ell\text{ prime}}{r_{\ell}(k-p^e)}, & \text{if }k>p^e, \\ 1, & \text{if }k=p^e.\end{cases}}
\]
Technically speaking, the algorithm only works when $m$ is the $n$-th partition number matrix for some $n\leqslant25012$, but the argument range can be extended provided that the list {\tt PrimesExtended} is extended accordingly.
\end{itemize}

\subsection{{\tt PartitionNumberMatrix}}\label{subsec2P3}

\begin{itemize}
\item Syntax: {\tt PartitionNumberMatrix(n)}, $n\in\IN^+$.
\item Output: the $n$-th partition number matrix.
\item Comments: Simply applies {\tt NextPartitionNumberMatrix} iteratively.
\end{itemize}

\subsection{{\tt NrPartitionsIntoPairwiseCoprimeNontrivialPrimePowers}}\label{subsec2P4}

\begin{itemize}
\item Syntax: {\tt NrPartitionsIntoPairwiseCoprimeNontrivialPrimePowers(n)}, $n\in\IN^+$.
\item Output: the number $r(n)$ of (un)ordered integer partitions of $n$ into pairwise coprime prime powers greater than $1$.
\item Comments: Computes the $n$-th partition number matrix and sums over the entries in its $n$-th row.
\end{itemize}

\subsection{{\tt NrElementOrdersSym}}\label{subsec2P5}

\begin{itemize}
\item Syntax: {\tt NrElementOrdersSym(n)}, $n\in\IN^+$.
\item Output: the number $\omicron(\Sym(n))$ of element orders in $\Sym(n)$.
\item Comments: This uses the fact that $\omicron(\Sym(n))=1+\sum_{k=1}^n{r(k)}$.
\end{itemize}

\subsection{{\tt ConstantsInOmicronSymBoundList}}\label{subsec2P6}

\begin{itemize}
\item Syntax: {\tt ConstantsInOmicronSymBoundList(n)}, $n\in\IN^+$.
\item Output: the list of the numbers $c_k:=\frac{\log{\omicron(\Sym(k))}}{\sqrt{k}}$ for $k=1,\ldots,n$.
\item Comments: By definition, $\omicron(\Sym(k))=\exp(c_k\sqrt{k})$. The authors conjecture that the maximum value of $c_k$ for all $k\in\IN^+$ is attained at $k=66$, and using this function, the authors verified this conjecture for $k\leqslant25000$.
\end{itemize}

\section{Functions for dealing with simple Lie type groups}\label{sec3}

\subsection{General-purpose functions}\label{subsec3P1}

\subsubsection{{\tt nrDivisors}}\label{subsubsec3P1P1}

\begin{itemize}
\item Syntax: {\tt nrDivisors(n)}, $n\in\IN^+$.
\item Output: the number of positive integer divisors of $n$.
\item Comments: Rather than just calling {\tt Length(DivisorsInt(n))}, which would involve computing the (possibly very long) list of all divisors of $n$, it calls {\tt FactorsInt(n)} to compute the factorisation of $n$ into pairwise coprime prime powers, $n=p_1^{k_1}\cdots p_r^{k_r}$, and then returns the product $(k_1+1)(k_2+1)\cdots(k_r+1)$.
\end{itemize}

\subsubsection{{\tt lcm}}\label{subsubsec3P1P2}

\begin{itemize}
\item Syntax: {\tt lcm(l)}, $l$ a list of positive integers.
\item Output: the least common multiple of the elements of {\tt l}.
\item Comments: Unlike GAP's built-in function {\tt Lcm}, our function {\tt lcm} returns $1$ (not an error message) when applied to the empty list.
\end{itemize}

\subsection{Functions for groups of type \texorpdfstring{$A_d$}{Ad}}\label{subsec3P2}

\subsubsection{{\tt OrderA}}\label{subsubsec3P2P1}

\begin{itemize}
\item Syntax: {\tt OrderA(d,q)}, $d\in\IN^+$, $q$ a prime power.
\item Output: the group order $|A_d(q)|=|\PSL_{d+1}(q)|$.
\end{itemize}

\subsubsection{{\tt LogLogOrderA}}\label{subsubsec3P2P2}

\begin{itemize}
\item Syntax: {\tt LogLogOrderA(d,q)}, $d\in\IN^+$, $q$ a prime power.
\item Output: the real number $\log\log{|A_d(q)|}$ as a decimal floating-point number.
\item Comments: This does \emph{not} just call {\tt Log(Log(Float(OrderA(d,q))))}, because applying {\tt Float} to too large integers does not result in an appropriate floating-point representation of those integers, but in the output {\tt inf}, which is treated like infinity.
\end{itemize}

\subsubsection{{\tt OrderOutA}}\label{subsubsec3P2P3}

\begin{itemize}
\item Syntax: {\tt OrderOutA(d,q)}, $d\in\IN^+$, $q$ a prime power.
\item Output: the outer automorphism group order $|\Out(A_d(q))|$.
\end{itemize}

\subsubsection{{\tt NrConjugacyClassesA}}\label{subsubsec3P2P4}

\begin{itemize}
\item Syntax: {\tt NrConjugacyClassesA(d,q)}, $d\in\IN^+$, $q$ a prime power.
\item Output: the conjugacy class number $\k(A_d(q))$.
\item Comments: Simply calls {\tt NrConjugacyClassesPSL(d+1,q)}.
\end{itemize}

\subsubsection{{\tt NrAutOrbitsBoundA}}\label{subsubsec3P2P5}

\begin{itemize}
\item Syntax: {\tt NrAutOrbitsBoundA(d,q,l)}, $d\in\IN^+$, $q$ a prime power, $l\in\IN^+$ ($l$ is interpreted as a \enquote{quality level}).
\item Output: If $l\not=2$, it returns the message \enquote{This quality level is not available. Please set the quality level to 2.} as a string. If $l=2$, it returns the lower bound $\underline{\omega}^{(2)}(A_d(q)):=\lceil\frac{\k(A_d(q))}{|\Out(A_d(q))|}\rceil$ on $\omega(A_d(q))$.
\end{itemize}

\subsubsection{{\tt epsilonOmegaBoundA}}\label{subsubsec3P2P6}

\begin{itemize}
\item Syntax: {\tt epsilonOmegaBoundA(d,q,l)}, $d\in\IN^+$, $q$ a prime power, $l\in\IN^+$ ($l$ is interpreted as a \enquote{quality level}).
\item Output: If $l\not=2$, it returns the message \enquote{This quality level is not available. Please set the quality level to 2.} as a string. If $l=2$, it returns the lower bound $\underline{\epsilon_{\omega}}^{(2)}(A_d(q)):=\frac{\log\log{\underline{\omega}^{(2)}(A_d(q))}}{\log\log{|A_d(q)|}}$ on $\epsilon_{\omega}(A_d(q))$.
\end{itemize}

\subsubsection{{\tt CoxeterNoA}}\label{subsubsec3P2P7}

\begin{itemize}
\item Syntax: {\tt CoxeterNoA(d)}, $d\in\IN^+$.
\item Output: the Coxeter number $h(A_d)=d+1$.
\end{itemize}

\subsubsection{{\tt SemisimpleElementOrdersA}}\label{subsubsec3P2P8}

\begin{itemize}
\item Syntax: {\tt SemisimpleElementOrdersA(d,q)}, $d\in\IN^+$, $q=p^f$ a prime power.
\item Output: the set of semisimple (viz., not divisible by $p$) element orders in $A_d(q)$.
\item Comments: It joins the sets of element orders of certain maximal tori $T$ of $A_d(q)$, using \cite[Theorem 2.1]{BG07a} to compute the group exponent $\Exp(T)$.
\end{itemize}

\subsubsection{{\tt NrSemisimpleElementOrdersA}}\label{subsubsec3P2P9}

\begin{itemize}
\item Syntax: {\tt NrSemisimpleElementOrdersA(d,q)}, $d\in\IN^+$, $q=p^f$ a prime power.
\item Output: the number $\omicron_{\ssrm}(A_d(q))=:\overline{\omicron_{\ssrm}}^{(2)}(A_d(q))$ of semisimple (viz., not divisible by $p$) element orders in $A_d(q)$.
\item Comments: Simply calls {\tt Length(SemisimpleElementOrdersA(d,q))}.
\end{itemize}

\subsubsection{{\tt NrSemisimpleElementOrdersBoundA}}\label{subsubsec3P2P10}

\begin{itemize}
\item Syntax: {\tt NrSemisimpleElementOrdersBoundA(d,q)}, $d\in\IN^+$, $q=p^f$ a prime power.
\item Output: an upper bound $\overline{\omicron_{\ssrm}}^{(1)}(A_d(q))$ on the number of semisimple (viz., not divisible by $p$) element orders in $A_d(q)$.
\item Comments: Like {\tt SemisimpleElementOrdersA(d,q)}, it loops over certain maximal tori of $A_d(q)$, but it adds the numbers of element orders in the tori that it loops over, rather than joining the corresponding sets. This makes it faster than {\tt SemisimpleElementOrdersA(d,q)}.
\end{itemize}

\subsubsection{{\tt SemisimpleElementOrdersGL}}\label{subsubsec3P2P11}

\begin{itemize}
\item Syntax: {\tt SemisimpleElementOrdersGL(n,q)}, $n\in\IN^+$, $q=p^f$ a prime power.
\item Output: the set $\Ord_{\ssrm}(\GL_n(q))$ of semisimple (viz., not divisible by $p$) element orders in $\GL_n(q)$.
\item Comments: Similar to {\tt SemisimpleElementOrdersA(d,q)}, but simpler (because the formula for $\Exp(T)$ when $T$ is a maximal torus of $\GL_n(q)$ is simpler than for maximal tori of $A_d(q)$).
\end{itemize}

\subsubsection{{\tt NrElementOrdersSharplyDivisibleByPToEBoundA}}\label{subsubsec3P2P12}

\begin{itemize}
\item Syntax: {\tt NrElementOrdersSharplyDivisibleByPToEBoundA(d,q,e)}, $d\in\IN^+$, $q=p^f$ a prime power, $e\in\{0,1,\ldots,1+\lfloor\log_p(d)\rfloor\}$.
\item Output: an upper bound $\overline{\omicron_{\ssrm}^{(e)}}(A_d(q))$ on the number $\omicron_{\ssrm}^{(e)}(A_d(q))$ of element orders $o$ in $A_d(q)$ such that the largest power of $p$ dividing $o$ is $p^e$.
\item Comments: For an explanation, see \cite[Case (5) of the proof of Theorem 1.1.3(5) in Subsection 3.3]{BGP19a}.
\end{itemize}

\subsubsection{{\tt NrElementOrdersBoundA}}\label{subsubsec3P2P13}

\begin{itemize}
\item Syntax: {\tt NrElementOrdersBoundA(d,q,l)}, $d\in\IN^+$, $q$ a prime power, $l\in\IN^+$ ($l$ is interpreted as a \enquote{quality level}).
\item Output: If $l\notin\{1,2,3\}$, it returns the message \enquote{This quality level is not available. Please set the quality level to 1, 2 or 3.} as a string. Otherwise, it returns a certain upper bound $\overline{\omicron}^{(l)}(A_d(q))$ on the number of element orders $\omicron(A_d(q))$. The larger $l$ is, the better is that bound, but also the more costly it is to compute.
\item Comments: Write $q=p^f$. For $l=3$, this function returns the sum of the numbers returned by {\tt NrElementOrdersSharplyDivisibleByPToEBoundA(d,q,e)} where $e=0,1,\ldots,1+\lfloor\log_p(d)\rfloor$. For $l\in\{1,2\}$, it returns $\overline{\omicron_{\ssrm}}^{(l)}(A_d(q))\cdot(1+\lceil\log_p(d+1)\rceil)$, the product of an upper bound on the number of semisimple element orders $\omicron_{\ssrm}(A_d(q))$ and the exact number of element orders in $A_d(q)$ that are powers of $p$.
\end{itemize}

\subsubsection{{\tt epsilonQBoundA}}\label{subsubsec3P2P14}

\begin{itemize}
\item Syntax: {\tt epsilonQBoundA(d,q,l1,l2)}, $d\in\IN^+$, $q$ a prime power, $l_1,l_2\in\IN^+$ (the $l_i$ are interpreted as \enquote{quality levels}).
\item Output: If $(l_1,l_2)\notin\{(2,1),(2,2),(2,3)\}$, then it returns the message \enquote{This combination of quality levels is not available. Please set the quality levels to $(2,1)$, $(2,2)$ or $(2,3)$.} as a string. Otherwise, it returns the lower bound
\[
\underline{\epsilon_{\q}}^{(l_1,l_2)}(A_d(q)):=\frac{\log\log{(\frac{\underline{\omega}^{(l_1)}(A_d(q))}{\overline{\omicron}^{(l_2)}(A_d(q)}+3)}}{\log\log{|A_d(q)|}}
\]
on $\epsilon_{\q}(A_d(q))$ as a decimal floating-point number.
\end{itemize}

\subsubsection{{\tt epsilonQBoundAd2}}\label{subsubsec3P2P15}

\begin{itemize}
\item Syntax: {\tt epsilonQBoundAd2(d)}, $d\in\IN^+$ (the intended range of use is $d\in\{54,\ldots,90\}$).
\item Output: the real number
\[
\frac{\log\log{(\frac{2^{d-1}}{\omicron_{\ssrm}(A_{90}(2))\cdot(1+\lceil\log_2(d+1)\rceil)}+3)}}{\log\log{|A_d(2)|}},
\]
as a decimal floating-point number.
\item Comments: If $d\leq90$, the output is a lower bound on $\epsilon_{\q}(A_d(2))$. It is used in \cite[Case (3) of the proof of Theorem 1.1.3(5) in Subsection 3.3]{BGP19a}.
\end{itemize}

\subsection{Functions for groups of type \texorpdfstring{${^2}A_d$}{2Ad}}\label{subsec3P3}

\subsubsection{{\tt Order2A}}\label{subsubsec3P3P1}

\begin{itemize}
\item Syntax: {\tt Order2A(d,Q)}, $d\in\IN^+$, $Q=q^2$ for a prime power $q$.
\item Output: the group order $|{^2}A_d(Q)|=|\PSU_{d+1}(q)|$.
\end{itemize}

\subsubsection{{\tt LogLogOrder2A}}\label{subsubsec3P3P2}

\begin{itemize}
\item Syntax: {\tt LogLogOrder2A(d,Q)}, $d\in\IN^+$, $Q=q^2$ for a prime power $q$.
\item Output: the real number $\log\log{|{^2}A_d(Q)|}$ as a decimal floating-point number.
\item Comments: This does \emph{not} just call {\tt Log(Log(Float(Order2A(d,Q))))}, because applying {\tt Float} to too large integers does not result in an appropriate floating-point representation of those integers, but in the output {\tt inf}, which is treated like infinity.
\end{itemize}

\subsubsection{{\tt OrderOut2A}}\label{subsubsec3P3P3}

\begin{itemize}
\item Syntax: {\tt OrderOut2A(d,Q)}, $d\in\IN^+$, $Q=q^2$ for a prime power $q$.
\item Output: the outer automorphism group order $|\Out({^2}A_d(Q))|$.
\end{itemize}

\subsubsection{{\tt NrConjugacyClasses2A}}\label{subsubsec3P3P4}

\begin{itemize}
\item Syntax: {\tt NrConjugacyClasses2A(d,Q)}, $d\in\IN^+$, $Q=q^2$ for a prime power $q$.
\item Output: the conjugacy class number $\k({^2}A_d(Q))$.
\item Comments: Simply calls {\tt NrConjugacyClassesPSU(d+1,q)}.
\end{itemize}

\subsubsection{{\tt NrAutOrbitsBound2A}}\label{subsubsec3P3P5}

\begin{itemize}
\item Syntax: {\tt NrAutOrbitsBound2A(d,Q,l)}, $d\in\IN^+$, $Q=q^2$ for a prime power $q$, $l\in\IN^+$ ($l$ is interpreted as a \enquote{quality level}).
\item Output: If $l\not=2$, it returns the message \enquote{This quality level is not available. Please set the quality level to 2.} as a string. If $l=2$, it returns the lower bound $\underline{\omega}^{(2)}({^2}A_d(Q)):=\lceil\frac{\k({^2}A_d(Q))}{|\Out({^2}A_d(Q))|}\rceil$ on $\omega({^2}A_d(Q))$.
\end{itemize}

\subsubsection{{\tt epsilonOmegaBound2A}}\label{subsubsec3P3P6}

\begin{itemize}
\item Syntax: {\tt epsilonOmegaBound2A(d,Q,l)}, $d\in\IN^+$, $Q=q^2$ for a prime power $q$, $l\in\IN^+$ ($l$ is interpreted as a \enquote{quality level}).
\item Output: If $l\not=2$, it returns the message \enquote{This quality level is not available. Please set the quality level to 2.} as a string. If $l=2$, it returns the lower bound $\underline{\epsilon_{\omega}}^{(2)}({^2}A_d(Q)):=\frac{\log\log{\underline{\omega}^{(2)}({^2}A_d(Q))}}{\log\log{|{^2}A_d(Q)|}}$ on $\epsilon_{\omega}({^2}A_d(Q))$.
\end{itemize}

\subsubsection{{\tt SemisimpleElementOrders2A}}\label{subsubsec3P3P7}

\begin{itemize}
\item Syntax: {\tt SemisimpleElementOrders2A(d,Q)}, $d\in\IN^+$, $Q=p^{2f}$ a prime power square.
\item Output: the set of semisimple (viz., not divisible by $p$) element orders in ${^2}A_d(Q)$.
\item Comments: It joins the sets of element orders of certain maximal tori $T$ of ${^2}A_d(Q)$, using \cite[Theorem 2.2]{BG07a} to compute the group exponent $\Exp(T)$.
\end{itemize}

\subsubsection{{\tt NrSemisimpleElementOrders2A}}\label{subsubsec3P3P8}

\begin{itemize}
\item Syntax: {\tt NrSemisimpleElementOrders2A(d,Q)}, $d\in\IN^+$, $Q=p^{2f}$ a prime power square.
\item Output: the number $\omicron_{\ssrm}({^2}A_d(Q))=:\overline{\omicron_{\ssrm}}^{(2)}({^2}A_d(Q))$ of semisimple (viz., not divisible by $p$) element orders in ${^2}A_d(Q)$.
\item Comments: Simply calls {\tt Length(SemisimpleElementOrders2A(d,Q))}.
\end{itemize}

\subsubsection{{\tt NrSemisimpleElementOrdersBound2A}}\label{subsubsec3P3P9}

\begin{itemize}
\item Syntax: {\tt NrSemisimpleElementOrdersBound2A(d,Q)}, $d\in\IN^+$, $Q=p^{2f}$ a prime power square.
\item Output: an upper bound $\overline{\omicron_{\ssrm}}^{(1)}({^2}A_d(Q))$ on the number of semisimple (viz., not divisible by $p$) element orders in ${^2}A_d(Q)$.
\item Comments: Like {\tt SemisimpleElementOrders2A(d,Q)}, it loops over certain maximal tori of ${^2}A_d(Q)$, but it adds the numbers of element orders in the tori that it loops over, rather than joining the corresponding sets. This makes it faster than {\tt SemisimpleElementOrders2A(d,Q)}.
\end{itemize}

\subsubsection{{\tt SemisimpleElementOrdersGU}}\label{subsubsec3P3P10}

\begin{itemize}
\item Syntax: {\tt SemisimpleElementOrdersGU(n,q)}, $n\in\IN^+$, $q=p^f$ a prime power.
\item Output: the set $\Ord_{\ssrm}(\GU_n(q))$ of semisimple (viz., not divisible by $p$) element orders in $\GU_n(q)$.
\item Comments: Similar to {\tt SemisimpleElementOrders2A(d,Q)}, but simpler (because the formula for $\Exp(T)$ when $T$ is a maximal torus of $\GU_n(q)$ is simpler than for maximal tori of ${^2}A_d(Q)$).
\end{itemize}

\subsubsection{{\tt NrElementOrdersSharplyDivisibleByPToEBound2A}}\label{subsubsec3P3P11}

\begin{itemize}
\item Syntax: {\tt NrElementOrdersSharplyDivisibleByPToEBound2A(d,Q,e)}, $d\in\IN^+$, $Q=q^2$ for a prime power $q$, $e\in\{0,1,\ldots,1+\lfloor\log_p(d)\rfloor\}$.
\item Output: an upper bound $\overline{\omicron_{\ssrm}^{(e)}}({^2}A_d(Q))$ on the number $\omicron_{\ssrm}^{(e)}({^2}A_d(Q))$ of element orders $o$ in ${^2}A_d(Q)$ such that the largest power of $p$ dividing $o$ is $p^e$.
\item Comments: For an explanation, see \cite[Case (5) of the proof of Theorem 1.1.3(5) in Subsection 3.3]{BGP19a}.
\end{itemize}

\subsubsection{{\tt NrElementOrdersBound2A}}\label{subsubsec3P3P12}

\begin{itemize}
\item Syntax: {\tt NrElementOrdersBound2A(d,Q,l)}, $d\in\IN^+$, $Q=q^2$ for a prime power $q$, $l\in\IN^+$ ($l$ is interpreted as a \enquote{quality level}).
\item Output: If $l\notin\{1,2,3\}$, it returns the message \enquote{This quality level is not available. Please set the quality level to 1, 2 or 3.} as a string. Otherwise, it returns a certain upper bound $\overline{\omicron}^{(l)}({^2}A_d(Q))$ on the number of element orders $\omicron({^2}A_d(Q))$. The larger $l$ is, the better is that bound, but also the more costly it is to compute.
\item Comments: Write $Q=p^{2f}$. For $l=3$, this function returns the sum of the numbers returned by {\tt NrElementOrdersSharplyDivisibleByPToEBound2A(d,Q,e)} where $e=0,1,\ldots,1+\lfloor\log_p(d)\rfloor$. For $l\in\{1,2\}$, it returns $\overline{\omicron_{\ssrm}}^{(l)}({^2}A_d(Q))\cdot(1+\lceil\log_p(d+1)\rceil)$, the product of an upper bound on the number of semisimple element orders $\omicron_{\ssrm}({^2}A_d(Q))$ and the exact number of element orders in ${^2}A_d(Q)$ that are powers of $p$.
\end{itemize}

\subsubsection{{\tt epsilonQBound2A}}\label{subsubsec3P3P13}

\begin{itemize}
\item Syntax: {\tt epsilonQBoundA(d,Q,l1,l2)}, $d\in\IN^+$, $Q=q^2$ for a prime power $q$, $l_1,l_2\in\IN^+$ (the $l_i$ are interpreted as \enquote{quality levels}).
\item Output: If $(l_1,l_2)\notin\{(2,1),(2,2),(2,3)\}$, then it returns the message \enquote{This combination of quality levels is not available. Please set the quality levels to $(2,1)$, $(2,2)$ or $(2,3)$.} as a string. Otherwise, it returns the lower bound
\[
\underline{\epsilon_{\q}}^{(l_1,l_2)}({^2}A_d(Q)):=\frac{\log\log{(\frac{\underline{\omega}^{(l_1)}({^2}A_d(Q))}{\overline{\omicron}^{(l_2)}({^2}A_d(Q)}+3)}}{\log\log{|{^2}A_d(Q)|}}
\]
on $\epsilon_{\q}({^2}A_d(Q))$ as a decimal floating-point number.
\end{itemize}

\subsubsection{{\tt epsilonQBound2Ad4}}\label{subsubsec3P3P14}

\begin{itemize}
\item Syntax: {\tt epsilonQBound2Ad4(d)}, $d\in\IN^+$ (the intended range of use is $d\in\{54,\ldots,90\}$).
\item Output: the real number
\[
\frac{\log\log{(\frac{2^d/(2\cdot\gcd(d+1,3)^2)}{\omicron_{\ssrm}({^2}A_{90}(4))\cdot(1+\lceil\log_2(d+1)\rceil)}+3)}}{\log\log{|{^2}A_d(4)|}},
\]
as a decimal floating-point number.
\item Comments: If $d\leq90$, the output is a lower bound on $\epsilon_{\q}({^2}A_d(4))$. It is used in \cite[Case (3) of the proof of Theorem 1.1.3(5) in Subsection 3.3]{BGP19a}.
\end{itemize}

\subsection{Functions for groups of type \texorpdfstring{$B_d$}{Bd}}\label{subsec3P4}

\subsubsection{{\tt OrderB}}\label{subsubsec3P4P1}

\begin{itemize}
\item Syntax: {\tt OrderB(d,q)}, $d\in\IN^+$, $q$ a prime power.
\item Output: the group order $|B_d(q)|=|\POmega_{2d+1}(q)|$.
\end{itemize}

\subsubsection{{\tt LogLogOrderB}}\label{subsubsec3P4P2}

\begin{itemize}
\item Syntax: {\tt LogLogOrderB(d,q)}, $d\in\IN^+$, $q$ a prime power.
\item Output: the real number $\log\log{|B_d(q)|}$ as a decimal floating-point number.
\item Comments: This does \emph{not} just call {\tt Log(Log(Float(OrderB(d,q))))}, because applying {\tt Float} to too large integers does not result in an appropriate floating-point representation of those integers, but in the output {\tt inf}, which is treated like infinity.
\end{itemize}

\subsubsection{{\tt OrderOutB}}\label{subsubsec3P4P3}

\begin{itemize}
\item Syntax: {\tt OrderOutB(d,q)}, $d\in\IN^+$, $q$ a prime power.
\item Output: the outer automorphism group order $|\Out(B_d(q))|$.
\end{itemize}

\subsubsection{{\tt NrConjugacyClassesB}}\label{subsubsec3P4P4}

\begin{itemize}
\item Syntax: {\tt NrConjugacyClassesB(d,q)}, $d\in\IN^+$, $q$ a prime power.
\item Output: the conjugacy class number $\k(B_d(q))=:\underline{\k}^{(2)}(B_d(q))$.
\item Comments: Uses \cite[Theorems 3.13(1) and 3.19(1)]{FG12a}.
\end{itemize}

\subsubsection{{\tt NrConjugacyClassesBoundB}}\label{subsec3P4P5}

\begin{itemize}
\item Syntax: {\tt NrConjugacyClassesBoundB(d,q)}, $d\in\IN^+$, $q$ a prime power.
\item Output: the lower bound $\underline{\k}^{(1)}(B_d(q)):=\lceil\frac{q^d}{\gcd(2,q-1)}\rceil$ on $\k(B_d(q))$.
\item Comments: To see that this is indeed a lower bound on $\k(B_d(q))$, use \cite[Theorem 1.1(1)]{FG12a}.
\end{itemize}

\subsubsection{{\tt NrAutOrbitsBoundB}}\label{subsubsec3P4P6}

\begin{itemize}
\item Syntax: {\tt NrAutOrbitsBoundB(d,q,l)}, $d\in\IN^+$, $q$ a prime power, $l\in\IN^+$ ($l$ is interpreted as a \enquote{quality level}).
\item Output: If $l\notin\{1,2\}$, it returns the message \enquote{This quality level is not available. Please set the quality level to 1 or 2.} as a string. Otherwise, it returns the lower bound $\underline{\omega}^{(l)}(B_d(q)):=\lceil\frac{\underline{\k}^{(l)}(B_d(q))}{|\Out(B_d(q))|}\rceil$ on $\omega(B_d(q))$.
\end{itemize}

\subsubsection{{\tt epsilonOmegaBoundB}}\label{subsubsec3P4P7}

\begin{itemize}
\item Syntax: {\tt epsilonOmegaBoundB(d,q,l)}, $d\in\IN^+$, $q$ a prime power, $l\in\IN^+$ ($l$ is interpreted as a \enquote{quality level}).
\item Output: If $l\notin\{1,2\}$, it returns the message \enquote{This quality level is not available. Please set the quality level to 1 or 2.} as a string. Otherwise, it returns the lower bound $\underline{\epsilon_{\omega}}^{(l)}(B_d(q)):=\frac{\log\log{\underline{\omega}^{(l)}(B_d(q))}}{\log\log{|B_d(q)|}}$ on $\epsilon_{\omega}(B_d(q))$.
\end{itemize}

\subsubsection{{\tt CoxeterNoB}}\label{subsubsec3P4P8}

\begin{itemize}
\item Syntax: {\tt CoxeterNoB(d)}, $d\in\IN^+$.
\item Output: the Coxeter number $h(B_d)=2d$.
\end{itemize}

\subsubsection{{\tt SemisimpleElementOrdersB}}\label{subsubsec3P4P9}

\begin{itemize}
\item Syntax: {\tt SemisimpleElementOrdersB(d,q)}, $d\in\IN^+$, $q=p^f$ a prime power.
\item Output: the set of semisimple (viz., not divisible by $p$) element orders in $B_d(q)$.
\item Comments: It joins the sets of element orders of certain maximal tori $T$ of $B_d(q)$, using \cite[Theorems 3 and 4]{BG07a} to compute the group exponent $\Exp(T)$.
\end{itemize}

\subsubsection{{\tt NrSemisimpleElementOrdersB}}\label{subsubsec3P4P10}

\begin{itemize}
\item Syntax: {\tt NrSemisimpleElementOrdersB(d,q)}, $d\in\IN^+$, $q=p^f$ a prime power.
\item Output: the number $\omicron_{\ssrm}(B_d(q))=:\overline{\omicron_{\ssrm}}^{(2)}(B_d(q))$ of semisimple (viz., not divisible by $p$) element orders in $B_d(q)$.
\item Comments: Simply calls {\tt Length(SemisimpleElementOrdersB(d,q))}.
\end{itemize}

\subsubsection{{\tt NrSemisimpleElementOrdersBoundB}}\label{subsubsec3P4P11}

\begin{itemize}
\item Syntax: {\tt NrSemisimpleElementOrdersBoundB(d,q)}, $d\in\IN^+$, $q=p^f$ a prime power.
\item Output: an upper bound $\overline{\omicron_{\ssrm}}^{(1)}(B_d(q))$ on the number of semisimple (viz., not divisible by $p$) element orders in $B_d(q)$.
\item Comments: Like {\tt SemisimpleElementOrdersB(d,q)}, it loops over certain maximal tori of $B_d(q)$, but it adds the numbers of element orders in the tori that it loops over, rather than joining the corresponding sets. This makes it faster than {\tt SemisimpleElementOrdersB(d,q)}.
\end{itemize}

\subsubsection{{\tt SemisimpleElementOrdersGOOddDim}}\label{subsubsec3P4P12}

\begin{itemize}
\item Syntax: {\tt SemisimpleElementOrdersGOOddDim(n,q)}, $n\in2\IN^++1$, $q=p^f$ a prime power.
\item Output: the set $\Ord_{\ssrm}(\GO_n(q))$ of semisimple (viz., not divisible by $p$) element orders in $\GO_n(q)$.
\item Comments: Similar to {\tt SemisimpleElementOrdersB(d,q)}, but simpler (because the formula for $\Exp(T)$ when $T$ is a maximal torus of $\GO_n(q)$ is simpler than for maximal tori of $B_d(q)$).
\end{itemize}

\subsubsection{{\tt NrElementOrdersBoundB}}\label{subsubsec3P4P13}

\begin{itemize}
\item Syntax: {\tt NrElementOrdersBoundB(d,q,l)}, $d\in\IN^+$, $q$ a prime power, $l\in\IN^+$ ($l$ is interpreted as a \enquote{quality level}).
\item Output: If $l\notin\{1,2\}$, it returns the message \enquote{This quality level is not available. Please set the quality level to 1 or 2.} as a string. Otherwise, it returns a certain upper bound $\overline{\omicron}^{(l)}(B_d(q))$ on the number of element orders $\omicron(B_d(q))$. The larger $l$ is, the better is that bound, but also the more costly it is to compute.
\item Comments: Write $q=p^f$. For $l\in\{1,2\}$, the output is $\overline{\omicron_{\ssrm}}^{(l)}(B_d(q))\cdot(1+\lceil\log_p(2d)\rceil)$, the product of an upper bound on $\omicron_{\ssrm}(B_d(q))$ and the exact number of element orders in $B_d(q)$ that are powers of $p$.
\end{itemize}

\subsubsection{{\tt epsilonQBoundB}}\label{subsubsec3P4P14}

\begin{itemize}
\item Syntax: {\tt epsilonQBoundB(d,q,l1,l2)}, $d\in\IN^+$, $q$ a prime power, $l_1,l_2\in\IN^+$ (the $l_i$ are interpreted as \enquote{quality levels}).
\item Output: If $(l_1,l_2)\notin\{(1,1),(1,2),(2,1),(2,2)\}$, then it returns the message \enquote{This combination of quality levels is not available. Please set the quality levels to $(1,1)$, $(1,2)$, $(2,1)$ or $(2,2)$.} as a string. Otherwise, it returns the lower bound
\[
\underline{\epsilon_{\q}}^{(l_1,l_2)}(B_d(q)):=\frac{\log\log{(\frac{\underline{\omega}^{(l_1)}(B_d(q))}{\overline{\omicron}^{(l_2)}(B_d(q)}+3)}}{\log\log{|B_d(q)|}}
\]
on $\epsilon_{\q}(B_d(q))$ as a decimal floating-point number.
\end{itemize}

\subsubsection{{\tt epsilonQBoundBd2}}\label{subsubsec3P4P15}

\begin{itemize}
\item Syntax: {\tt epsilonQBoundBd2(d)}, $d\in\IN^+$ (the intended range of use is $d\in\{54,\ldots,90\}$).
\item Output: the real number
\[
\frac{\log\log{(\frac{\underline{\omega}^{(2)}(B_d(2))}{\omicron_{\ssrm}(B_{90}(2))\cdot(1+\lceil\log_2(2d)\rceil)}+3)}}{\log\log{|B_d(2)|}},
\]
as a decimal floating-point number.
\item Comments: If $d\leq90$, the output is a lower bound on $\epsilon_{\q}(B_d(2))$. It is used in \cite[Case (3) of the proof of Theorem 1.1.3(5) in Subsection 3.3]{BGP19a}.
\end{itemize}

\subsection{Functions for groups of type \texorpdfstring{$C_d$}{Cd}}\label{subsec3P5}

\subsubsection{{\tt OrderC}}\label{subsubsec3P5P1}

\begin{itemize}
\item Syntax: {\tt OrderC(d,q)}, $d\in\IN^+$, $q$ a prime power.
\item Output: the group order $|C_d(q)|=|\PSp_{2d}(q)|$.
\end{itemize}

\subsubsection{{\tt LogLogOrderC}}\label{subsubsec3P5P2}

\begin{itemize}
\item Syntax: {\tt LogLogOrderC(d,q)}, $d\in\IN^+$, $q$ a prime power.
\item Output: the real number $\log\log{|C_d(q)|}$ as a decimal floating-point number.
\item Comments: This does \emph{not} just call {\tt Log(Log(Float(OrderC(d,q))))}, because applying {\tt Float} to too large integers does not result in an appropriate floating-point representation of those integers, but in the output {\tt inf}, which is treated like infinity.
\end{itemize}

\subsubsection{{\tt OrderOutC}}\label{subsubsec3P5P3}

\begin{itemize}
\item Syntax: {\tt OrderOutC(d,q)}, $d\in\IN^+$, $q$ a prime power.
\item Output: the outer automorphism group order $|\Out(C_d(q))|$.
\end{itemize}

\subsubsection{{\tt NrConjugacyClassesSp}}\label{subsubsec3P5P4}

\begin{itemize}
\item Syntax: {\tt NrConjugacyClassesSp(n,q)}, $n\in2\IN^+$, $q$ a prime power.
\item Output: the conjugacy class number $\k(\Sp_n(q))$.
\item Comments: Uses \cite[Case (B), statement (iii), p.~36]{Wal63a} \cite[Theorem 3.13(1)]{FG12a}.
\end{itemize}

\subsubsection{{\tt NrConjugacyClassesBoundC}}\label{subsec3P5P5}

\begin{itemize}
\item Syntax: {\tt NrConjugacyClassesBoundC(d,q,l)}, $d\in\IN^+$, $q$ a prime power, $l\in\IN^+$ ($l$ is interpreted as a \enquote{quality level}).
\item Output: If $l\notin\{1,2\}$, it returns the message \enquote{This quality level is not available. Please set the quality level to 1 or 2.} as a string. If $l=2$, it returns the lower bound $\underline{\k}^{(2)}(C_d(q)):=\lceil\frac{\k(\Sp_{2d}(q))}{\gcd(2,q-1)}\rceil$ on $\k(C_d(q))$. If $l=1$, it returns the lower bound $\underline{\k}^{(1)}(C_d(q)):=\lceil\frac{q^d}{\gcd(2,q-1)}\rceil$ on $\k(C_d(q))$.
\item Comments: To see that $\k(\Inndiag(C_d(q)))\geqslant q^d$ (and thus that the output for $l=1$ is indeed a lower bound on $\k(C_d(q))$), use \cite[Theorem 1.1(1)]{FG12a}.
\end{itemize}

\subsubsection{{\tt NrAutOrbitsBoundC}}\label{subsubsec3P5P6}

\begin{itemize}
\item Syntax: {\tt NrAutOrbitsBoundC(d,q,l)}, $d\in\IN^+$, $q$ a prime power, $l\in\IN^+$ ($l$ is interpreted as a \enquote{quality level}).
\item Output: If $l\notin\{1,2\}$, it returns the message \enquote{This quality level is not available. Please set the quality level to 1 or 2.} as a string. Otherwise, it returns the lower bound $\underline{\omega}^{(l)}(C_d(q)):=\lceil\frac{\underline{\k}^{(l)}(C_d(q))}{|\Out(C_d(q))|}\rceil$ on $\omega(C_d(q))$.
\end{itemize}

\subsubsection{{\tt epsilonOmegaBoundC}}\label{subsubsec3P5P7}

\begin{itemize}
\item Syntax: {\tt epsilonOmegaBoundC(d,q,l)}, $d\in\IN^+$, $q$ a prime power, $l\in\IN^+$ ($l$ is interpreted as a \enquote{quality level}).
\item Output: If $l\notin\{1,2\}$, it returns the message \enquote{This quality level is not available. Please set the quality level to 1 or 2.} as a string. Otherwise, it returns the lower bound $\underline{\epsilon_{\omega}}^{(l)}(C_d(q)):=\frac{\log\log{\underline{\omega}^{(l)}(C_d(q))}}{\log\log{|C_d(q)|}}$ on $\epsilon_{\omega}(C_d(q))$.
\end{itemize}

\subsubsection{{\tt CoxeterNoC}}\label{subsubsec3P5P8}

\begin{itemize}
\item Syntax: {\tt CoxeterNoC(d)}, $d\in\IN^+$.
\item Output: the Coxeter number $h(C_d)=2d$.
\end{itemize}

\subsubsection{{\tt SemisimpleElementOrdersC}}\label{subsubsec3P5P9}

\begin{itemize}
\item Syntax: {\tt SemisimpleElementOrdersC(d,q)}, $d\in\IN^+$, $q=p^f$ a prime power.
\item Output: the set of semisimple (viz., not divisible by $p$) element orders in $C_d(q)$.
\item Comments: By \cite[Theorems 3 and 4]{BG07a}, the sets of semisimple element orders in $B_d(q)$ and $C_d(q)$ respectively are equal. Hence this function simply calls {\tt SemisimpleElementOrdersB(d,q)}.
\end{itemize}

\subsubsection{{\tt NrSemisimpleElementOrdersC}}\label{subsubsec3P5P10}

\begin{itemize}
\item Syntax: {\tt NrSemisimpleElementOrdersC(d,q)}, $d\in\IN^+$, $q=p^f$ a prime power.
\item Output: the number $\omicron_{\ssrm}(C_d(q))=:\overline{\omicron_{\ssrm}}^{(2)}(C_d(q))$ of semisimple (viz., not divisible by $p$) element orders in $C_d(q)$.
\item Comments: Simply calls {\tt Length(SemisimpleElementOrdersC(d,q))}.
\end{itemize}

\subsubsection{{\tt NrSemisimpleElementOrdersBoundC}}\label{subsubsec3P5P11}

\begin{itemize}
\item Syntax: {\tt NrSemisimpleElementOrdersBoundC(d,q)}, $d\in\IN^+$, $q=p^f$ a prime power.
\item Output: an upper bound $\overline{\omicron_{\ssrm}}^{(1)}(C_d(q))$ on the number of semisimple (viz., not divisible by $p$) element orders in $C_d(q)$.
\item Comments: Simply calls {\tt NrSemisimpleElementOrdersBoundB(d,q)} (cf.~also the comments in Subsection \ref{subsubsec3P5P9}).
\end{itemize}

\subsubsection{{\tt NrElementOrdersBoundC}}\label{subsubsec3P5P12}

\begin{itemize}
\item Syntax: {\tt NrElementOrdersBoundC(d,q,l)}, $d\in\IN^+$, $q$ a prime power, $l\in\IN^+$ ($l$ is interpreted as a \enquote{quality level}).
\item Output: If $l\notin\{1,2\}$, it returns the message \enquote{This quality level is not available. Please set the quality level to 1 or 2.} as a string. Otherwise, it returns a certain upper bound $\overline{\omicron}^{(l)}(C_d(q))$ on the number of element orders $\omicron(C_d(q))$. The larger $l$ is, the better is that bound, but also the more costly it is to compute.
\item Comments: Write $q=p^f$. For $l\in\{1,2\}$, the output is $\overline{\omicron_{\ssrm}}^{(l)}(C_d(q))\cdot(1+\lceil\log_p(2d)\rceil)$, the product of an upper bound on $\omicron_{\ssrm}(C_d(q))$ and the exact number of element orders in $C_d(q)$ that are powers of $p$.
\end{itemize}

\subsubsection{{\tt epsilonQBoundC}}\label{subsubsec3P5P13}

\begin{itemize}
\item Syntax: {\tt epsilonQBoundC(d,q,l1,l2)}, $d\in\IN^+$, $q$ a prime power, $l_1,l_2\in\IN^+$ (the $l_i$ are interpreted as \enquote{quality levels}).
\item Output: If $(l_1,l_2)\notin\{(1,1),(1,2),(2,1),(2,2)\}$, then it returns the message \enquote{This combination of quality levels is not available. Please set the quality levels to $(1,1)$, $(1,2)$, $(2,1)$ or $(2,2)$.} as a string. Otherwise, it returns the lower bound
\[
\underline{\epsilon_{\q}}^{(l_1,l_2)}(C_d(q)):=\frac{\log\log{(\frac{\underline{\omega}^{(l_1)}(C_d(q))}{\overline{\omicron}^{(l_2)}(C_d(q)}+3)}}{\log\log{|C_d(q)|}}
\]
on $\epsilon_{\q}(B_d(q))$ as a decimal floating-point number.
\end{itemize}

\subsubsection{{\tt epsilonQBoundCd2}}\label{subsubsec3P5P14}

\begin{itemize}
\item Syntax: {\tt epsilonQBoundCd2(d)}, $d\in\IN^+$ (the intended range of use is $d\in\{54,\ldots,90\}$).
\item Output: the real number
\[
\frac{\log\log{(\frac{\underline{\omega}^{(2)}(C_d(2))}{\omicron_{\ssrm}(C_{90}(2))\cdot(1+\lceil\log_2(2d)\rceil)}+3)}}{\log\log{|C_d(2)|}},
\]
as a decimal floating-point number.
\item Comments: If $d\leq90$, the output is a lower bound on $\epsilon_{\q}(C_d(2))$. It is used in \cite[Case (3) of the proof of Theorem 1.1.3(5) in Subsection 3.3]{BGP19a}.
\end{itemize}

\subsection{Functions for groups of type \texorpdfstring{$D_d$}{Dd}}\label{subsec3P6}

\subsubsection{{\tt OrderD}}\label{subsubsec3P6P1}

\begin{itemize}
\item Syntax: {\tt OrderD(d,q)}, $d\in\IN^+$, $d\geqslant2$, $q$ a prime power.
\item Output: the group order $|D_d(q)|=|\POmega^+_{2d}(q)|$.
\end{itemize}

\subsubsection{{\tt LogLogOrderD}}\label{subsubsec3P6P2}

\begin{itemize}
\item Syntax: {\tt LogLogOrderD(d,q)}, $d\in\IN^+$, $d\geqslant2$, $q$ a prime power.
\item Output: the real number $\log\log{|D_d(q)|}$ as a decimal floating-point number.
\item Comments: This does \emph{not} just call {\tt Log(Log(Float(OrderD(d,q))))}, because applying {\tt Float} to too large integers does not result in an appropriate floating-point representation of those integers, but in the output {\tt inf}, which is treated like infinity.
\end{itemize}

\subsubsection{{\tt OrderOutD}}\label{subsubsec3P6P3}

\begin{itemize}
\item Syntax: {\tt OrderOutD(d,q)}, $d\in\IN^+$, $d\geqslant2$, $q$ a prime power.
\item Output: the outer automorphism group order $|\Out(D_d(q))|$.
\end{itemize}

\subsubsection{{\tt NrConjugacyClassesSumOmega}}\label{subsubsec3P6P4}

\begin{itemize}
\item Syntax: {\tt NrConjugacyClassesSumOmega(n,q)}, $n\in2\IN^++2$, $q$ an \emph{even} prime power.
\item Output: the sum $\k(\Omega^+_n(q))+\k(\Omega^-_n(q))$.
\item Comments: See \cite[Theorem 3.22(1)]{FG12a}.
\end{itemize}

\subsubsection{{\tt NrConjugacyClassesDifferenceOmega}}\label{subsubsec3P6P5}

\begin{itemize}
\item Syntax: {\tt NrConjugacyClassesDifferenceOmega(n,q)}, $n\in2\IN^++2$, $q$ an \emph{even} prime power.
\item Output: the difference $\k(\Omega^+_n(q))-\k(\Omega^-_n(q))$.
\item Comments: See \cite[Theorem 3.22(2)]{FG12a}.
\end{itemize}

\subsubsection{{\tt NrConjugacyClassesSumSO}}\label{subsubec3P6P6}

\begin{itemize}
\item Syntax: {\tt NrConjugacyClassesSumSO(n,q)}, $n\in2\IN^++2$, $q$ an \emph{odd} prime power.
\item Output: the sum $\k(\SO^+_n(q))+\k(\Omega^-_n(q))$.
\item Comments: See \cite[Theorem 3.16(1)]{FG12a}.
\end{itemize}

\subsubsection{{\tt NrConjugacyClassesDifferenceSO}}\label{subsubec3P6P7}

\begin{itemize}
\item Syntax: {\tt NrConjugacyClassesDifferenceSO(n,q)}, $n\in2\IN^++2$, $q$ an \emph{odd} prime power.
\item Output: the sum $\k(\SO^+_n(q))-\k(\Omega^-_n(q))$.
\item Comments: See \cite[Theorem 3.16(2)]{FG12a}.
\end{itemize}

\subsubsection{{\tt NrConjugacyClassesSO}}\label{subsubec3P6P8}

\begin{itemize}
\item Syntax: {\tt NrConjugacyClassesSO(epsilon,n,q)}, $\epsilon\in\{\pm1\}$, $n\in2\IN^++2$, $q$ an \emph{odd} prime power.
\item Output: $\k(\SO^{\sign(\epsilon)}_n(q))$.
\end{itemize}

\subsubsection{{\tt NrConjugacyClassesOmega}}\label{subsubsec3P6P9}

\begin{itemize}
\item Syntax: {\tt NrConjugacyClassesOmega(epsilon,n,q)}, $\epsilon\in\{\pm1\}$, $n\in2\IN^++2$, $q$ a prime power.
\item Output: $\k(\Omega^{\sign(\epsilon)}_n(q))$.
\item Comments: For odd $q$, this uses \cite[Theorem 3.18(1) and the paragraph before Theorem 3.18]{FG12a}.
\end{itemize}

\subsubsection{{\tt NrConjugacyClassesBoundD}}\label{subsec3P6P10}

\begin{itemize}
\item Syntax: {\tt NrConjugacyClassesBoundD(d,q,l)}, $d\in\IN^+$, $d\geqslant2$, $q$ a prime power, $l\in\IN^+$ ($l$ is interpreted as a \enquote{quality level}).
\item Output: If $l\notin\{1,2\}$, it returns the message \enquote{This quality level is not available. Please set the quality level to 1 or 2.} as a string. If $l=2$, it returns the following lower bound $\underline{\k}^{(2)}(D_d(q))$ on $\k(D_d(q))$:
\[
\underline{\k}^{(2)}(D_d(q)):=\begin{cases}\k(\Omega^+_{2d}(q)), & \text{if }2\mid q, \\ \frac{1}{2}\k(\SO^+_{2d}(q)), & \text{if }2\nmid q,q\equiv3\Mod{4}\text{ and }2\nmid d, \\ \lceil\frac{1}{2}\k(\Omega^+_{2d}(q))\rceil, & \text{if }2\nmid q,\text{ and either}q\equiv1\Mod{4}\text{ or }2\mid d.\end{cases}
\]
If $l=1$, it returns the lower bound $\underline{\k}^{(1)}(D_d(q)):=\frac{q^d}{\gcd(2,q-1)^2}$ on $\k(D_d(q))$.
\end{itemize}

\subsubsection{{\tt NrAutOrbitsBoundD}}\label{subsubsec3P6P11}

\begin{itemize}
\item Syntax: {\tt NrAutOrbitsBoundD(d,q,l)}, $d\in\IN^+$, $d\geqslant2$, $q$ a prime power, $l\in\IN^+$ ($l$ is interpreted as a \enquote{quality level}).
\item Output: If $l\notin\{1,2\}$, it returns the message \enquote{This quality level is not available. Please set the quality level to 1 or 2.} as a string. Otherwise, it returns the lower bound $\underline{\omega}^{(l)}(D_d(q)):=\lceil\frac{\underline{\k}^{(l)}(D_d(q))}{|\Out(D_d(q))|}\rceil$ on $\omega(D_d(q))$.
\end{itemize}

\subsubsection{{\tt epsilonOmegaBoundD}}\label{subsubsec3P6P12}

\begin{itemize}
\item Syntax: {\tt epsilonOmegaBoundD(d,q,l)}, $d\in\IN^+$, $d\geqslant2$, $q$ a prime power, $l\in\IN^+$ ($l$ is interpreted as a \enquote{quality level}).
\item Output: If $l\notin\{1,2\}$, it returns the message \enquote{This quality level is not available. Please set the quality level to 1 or 2.} as a string. Otherwise, it returns the lower bound $\underline{\epsilon_{\omega}}^{(l)}(D_d(q)):=\frac{\log\log{\underline{\omega}^{(l)}(D_d(q))}}{\log\log{|D_d(q)|}}$ on $\epsilon_{\omega}(D_d(q))$.
\end{itemize}

\subsubsection{{\tt CoxeterNoD}}\label{subsubsec3P6P13}

\begin{itemize}
\item Syntax: {\tt CoxeterNoD(d)}, $d\in\IN^+$, $d\geqslant2$.
\item Output: the Coxeter number $h(D_d)=2d-2$.
\end{itemize}

\subsubsection{{\tt SemisimpleElementOrdersD}}\label{subsubsec3P614}

\begin{itemize}
\item Syntax: {\tt SemisimpleElementOrdersD(d,q)}, $d\in\IN^+$, $d\geqslant2$, $q=p^f$ a prime power.
\item Output: the set of semisimple (viz., not divisible by $p$) element orders in $D_d(q)$.
\item Comments: It joins the sets of element orders of certain maximal tori $T$ of $D_d(q)$, using \cite[Theorems 5 and 7]{BG07a} to compute the group exponent $\Exp(T)$.
\end{itemize}

\subsubsection{{\tt NrSemisimpleElementOrdersD}}\label{subsubsec3P6P15}

\begin{itemize}
\item Syntax: {\tt NrSemisimpleElementOrdersD(d,q)}, $d\in\IN^+$, $d\geqslant2$, $q=p^f$ a prime power.
\item Output: the number $\omicron_{\ssrm}(D_d(q))=:\overline{\omicron_{\ssrm}}^{(2)}(D_d(q))$ of semisimple (viz., not divisible by $p$) element orders in $D_d(q)$.
\item Comments: Simply calls {\tt Length(SemisimpleElementOrdersD(d,q))}.
\end{itemize}

\subsubsection{{\tt NrSemisimpleElementOrdersBoundD}}\label{subsubsec3P6P16}

\begin{itemize}
\item Syntax: {\tt NrSemisimpleElementOrdersBoundD(d,q)}, $d\in\IN^+$, $d\geqslant2$, $q=p^f$ a prime power.
\item Output: an upper bound $\overline{\omicron_{\ssrm}}^{(1)}(D_d(q))$ on the number of semisimple (viz., not divisible by $p$) element orders in $D_d(q)$.
\item Comments: Like {\tt SemisimpleElementOrdersD(d,q)}, it loops over certain maximal tori of $D_d(q)$, but it adds the numbers of element orders in the tori that it loops over, rather than joining the corresponding sets. This makes it faster than {\tt SemisimpleElementOrdersD(d,q)}.
\end{itemize}

\subsubsection{{\tt SemisimpleElementOrdersGOEvenDim}}\label{subsubsec3P6P17}

\begin{itemize}
\item Syntax: {\tt SemisimpleElementOrdersGOEvenDim(epsilon,n,q)}, $\epsilon\in\{\pm1\}$, $n\in2\IN^++2$, $q=p^f$ a prime power.
\item Output: the set $\Ord_{\ssrm}(\GO^{\sign(\epsilon)}_n(q))$ of semisimple (viz., not divisible by $p$) element orders in $\GO^{\sign(\epsilon)}_n(q)$.
\end{itemize}

\subsubsection{{\tt NrElementOrdersSharplyDivisibleByPToEBoundD}}\label{subsubsec3P6P18}

\begin{itemize}
\item Syntax: {\tt NrElementOrdersSharplyDivisibleByPToEBoundD(d,q,e)}, $d\in\IN^+$, $d\geqslant2$, $q=p^f$ a prime power, $e\in\{0,1,\ldots,1+\lceil\log_p(2d-3)\rceil\}$.
\item Output: an upper bound $\overline{\omicron_{\ssrm}^{(e)}}(D_d(q))$ on the number $\omicron_{\ssrm}^{(e)}(D_d(q))$ of element orders $o$ in $D_d(q)$ such that the largest power of $p$ dividing $o$ is $p^e$.
\item Comments: For an explanation, see \cite[Case (5) of the proof of Theorem 1.1.3(5) in Subsection 3.3]{BGP19a}.
\end{itemize}

\subsubsection{{\tt NrElementOrdersBoundD}}\label{subsubsec3P6P19}

\begin{itemize}
\item Syntax: {\tt NrElementOrdersBoundD(d,q,l)}, $d\in\IN^+$, $d\geqslant2$, $q$ a prime power, $l\in\IN^+$ ($l$ is interpreted as a \enquote{quality level}).
\item Output: If $l\notin\{1,2,3\}$, it returns the message \enquote{This quality level is not available. Please set the quality level to 1, 2 or 3.} as a string. Otherwise, it returns a certain upper bound $\overline{\omicron}^{(l)}(D_d(q))$ on the number of element orders $\omicron(D_d(q))$. The larger $l$ is, the better is that bound, but also the more costly it is to compute.
\item Comments: Write $q=p^f$. For $l=3$, this function returns the sum of the numbers returned by {\tt NrElementOrdersSharplyDivisibleByPToEBoundD(d,q,e)} where $e=0,1,\ldots,1+\lfloor\log_p(2d-3)\rfloor$. For $l\in\{1,2\}$, the output is $\overline{\omicron_{\ssrm}}^{(l)}(D_d(q))\cdot(1+\lceil\log_p(2d-2)\rceil)$, the product of an upper bound on $\omicron_{\ssrm}(D_d(q))$ and the exact number of element orders in $D_d(q)$ that are powers of $p$.
\end{itemize}

\subsubsection{{\tt epsilonQBoundD}}\label{subsubsec3P6P20}

\begin{itemize}
\item Syntax: {\tt epsilonQBoundD(d,q,l1,l2)}, $d\in\IN^+$, $d\geqslant2$, $q$ a prime power, $l_1,l_2\in\IN^+$ (the $l_i$ are interpreted as \enquote{quality levels}).
\item Output: If $(l_1,l_2)\notin\{1,2\}\times\{1,2,3\}$, then it returns the message \enquote{This combination of quality levels is not available. Please set the quality levels to $(1,1)$, $(1,2)$, $(1,3)$, $(2,1)$, $(2,2)$ or $(2,3)$.} as a string. Otherwise, it returns the lower bound
\[
\underline{\epsilon_{\q}}^{(l_1,l_2)}(D_d(q)):=\frac{\log\log{(\frac{\underline{\omega}^{(l_1)}(D_d(q))}{\overline{\omicron}^{(l_2)}(D_d(q)}+3)}}{\log\log{|D_d(q)|}}
\]
on $\epsilon_{\q}(D_d(q))$ as a decimal floating-point number.
\end{itemize}

\subsubsection{{\tt epsilonQBoundDd2}}\label{subsubsec3P6P21}

\begin{itemize}
\item Syntax: {\tt epsilonQBoundDd2(d)}, $d\in\IN^+$ (the intended range of use is $d\in\{54,\ldots,90\}$).
\item Output: the real number
\[
\frac{\log\log{(\frac{\underline{\omega}^{(2)}(D_d(2))}{\omicron_{\ssrm}(D_{90}(2))\cdot(1+\lceil\log_2(2d-2)\rceil)}+3)}}{\log\log{|D_d(2)|}},
\]
as a decimal floating-point number.
\item Comments: If $d\leq90$, the output is a lower bound on $\epsilon_{\q}(D_d(2))$. It is used in \cite[Case (3) of the proof of Theorem 1.1.3(5) in Subsection 3.3]{BGP19a}.
\end{itemize}

\subsection{Functions for groups of type \texorpdfstring{${^2}D_d$}{2Dd}}\label{subsec3P7}

\subsubsection{{\tt Order2D}}\label{subsubsec3P7P1}

\begin{itemize}
\item Syntax: {\tt Order2D(d,q)}, $d\in\IN^+$, $d\geqslant2$, $Q=q^2$ for a prime power $q$.
\item Output: the group order $|{^2}D_d(Q)|=|\POmega^-_{2d}(q)|$.
\end{itemize}

\subsubsection{{\tt LogLogOrder2D}}\label{subsubsec3P7P2}

\begin{itemize}
\item Syntax: {\tt LogLogOrder2D(d,Q)}, $d\in\IN^+$, $d\geqslant2$, $Q=q^2$ for a prime power $q$.
\item Output: the real number $\log\log{|{^2}D_d(Q)|}$ as a decimal floating-point number.
\item Comments: This does \emph{not} just call {\tt Log(Log(Float(Order2D(d,Q))))}, because applying {\tt Float} to too large integers does not result in an appropriate floating-point representation of those integers, but in the output {\tt inf}, which is treated like infinity.
\end{itemize}

\subsubsection{{\tt OrderOut2D}}\label{subsubsec3P7P3}

\begin{itemize}
\item Syntax: {\tt OrderOut2D(d,Q)}, $d\in\IN^+$, $d\geqslant2$, $Q=q^2$ a prime power.
\item Output: the outer automorphism group order $|\Out({^2}D_d(Q))|$.
\end{itemize}

\subsubsection{{\tt NrConjugacyClassesBound2D}}\label{subsec3P7P4}

\begin{itemize}
\item Syntax: {\tt NrConjugacyClassesBound2D(d,Q,l)}, $d\in\IN^+$, $d\geqslant2$, $Q=q^2$ for a prime power $q$, $l\in\IN^+$ ($l$ is interpreted as a \enquote{quality level}).
\item Output: If $l\notin\{1,2\}$, it returns the message \enquote{This quality level is not available. Please set the quality level to 1 or 2.} as a string. If $l=2$, it returns the following lower bound $\underline{\k}^{(2)}({^2}D_d(Q))$ on $\k({^2}D_d(Q))$:
\[
\underline{\k}^{(2)}({^2}D_d(Q)):=\begin{cases}\k(\Omega^-_{2d}(q)), & \text{if }2\mid q, \\ \frac{1}{2}\k(\SO^-_{2d}(q)), & \text{if }2\nmid q,\text{ and either }q\equiv1\Mod{4}\text{ or }2\mid d, \\ \lceil\frac{1}{2}\k(\Omega^-_{2d}(q))\rceil, & \text{if }2\nmid q,q\equiv3\Mod{4}\text{ and }2\nmid d.\end{cases}
\]
If $l=1$, it returns the lower bound $\underline{\k}^{(1)}(^{2}D_d(Q)):=\frac{q^d}{\gcd(2,q-1)^2}$ on $\k(^{2}D_d(Q))$.
\end{itemize}

\subsubsection{{\tt NrAutOrbitsBound2D}}\label{subsubsec3P7P5}

\begin{itemize}
\item Syntax: {\tt NrAutOrbitsBound2D(d,Q,l)}, $d\in\IN^+$, $d\geqslant2$, $Q=q^2$ for a prime power $q$, $l\in\IN^+$ ($l$ is interpreted as a \enquote{quality level}).
\item Output: If $l\notin\{1,2\}$, it returns the message \enquote{This quality level is not available. Please set the quality level to 1 or 2.} as a string. Otherwise, it returns the lower bound $\underline{\omega}^{(l)}({^2}D_d(Q)):=\lceil\frac{\underline{\k}^{(l)}({^2}D_d(Q))}{|\Out({^2}D_d(Q))|}\rceil$ on $\omega({^2}D_d(Q))$.
\end{itemize}

\subsubsection{{\tt epsilonOmegaBound2D}}\label{subsubsec3P7P6}

\begin{itemize}
\item Syntax: {\tt epsilonOmegaBound2D(d,Q,l)}, $d\in\IN^+$, $d\geqslant2$, $Q=q^2$ for a prime power $q$, $l\in\IN^+$ ($l$ is interpreted as a \enquote{quality level}).
\item Output: If $l\notin\{1,2\}$, it returns the message \enquote{This quality level is not available. Please set the quality level to 1 or 2.} as a string. Otherwise, it returns the lower bound $\underline{\epsilon_{\omega}}^{(l)}({^2}D_d(Q)):=\frac{\log\log{\underline{\omega}^{(l)}({^2}D_d(Q))}}{\log\log{|{^2}D_d(Q)|}}$ on $\epsilon_{\omega}({^2}D_d(Q))$.
\end{itemize}

\subsubsection{{\tt SemisimpleElementOrders2D}}\label{subsubsec3P7P7}

\begin{itemize}
\item Syntax: {\tt SemisimpleElementOrders2D(d,Q)}, $d\in\IN^+$, $d\geqslant2$, $Q=p^{2f}$ a prime power square.
\item Output: the set of semisimple (viz., not divisible by $p$) element orders in ${^2}D_d(Q)$.
\item Comments: It joins the sets of element orders of certain maximal tori $T$ of ${^2}D_d(Q)$, using \cite[Theorems 6 and 7]{BG07a} to compute the group exponent $\Exp(T)$.
\end{itemize}

\subsubsection{{\tt NrSemisimpleElementOrders2D}}\label{subsubsec3P7P8}

\begin{itemize}
\item Syntax: {\tt NrSemisimpleElementOrders2D(d,Q)}, $d\in\IN^+$, $d\geqslant2$, $Q=p^{2f}$ a prime power square.
\item Output: the number $\omicron_{\ssrm}({^2}D_d(Q))=:\overline{\omicron_{\ssrm}}^{(2)}({^2}D_d(Q))$ of semisimple (viz., not divisible by $p$) element orders in ${^2}D_d(Q)$.
\item Comments: Simply calls {\tt Length(SemisimpleElementOrders2D(d,Q))}.
\end{itemize}

\subsubsection{{\tt NrSemisimpleElementOrdersBound2D}}\label{subsubsec3P7P9}

\begin{itemize}
\item Syntax: {\tt NrSemisimpleElementOrdersBound2D(d,Q)}, $d\in\IN^+$, $d\geqslant2$, $Q=p^{2f}$ a prime power square.
\item Output: an upper bound $\overline{\omicron_{\ssrm}}^{(1)}({^2}D_d(Q))$ on the number of semisimple (viz., not divisible by $p$) element orders in ${^2}D_d(Q)$.
\item Comments: Like {\tt SemisimpleElementOrders2D(d,Q)}, it loops over certain maximal tori of ${^2}D_d(Q)$, but it adds the numbers of element orders in the tori that it loops over, rather than joining the corresponding sets. This makes it faster than {\tt SemisimpleElementOrders2D(d,Q)}.
\end{itemize}

\subsubsection{{\tt NrElementOrdersSharplyDivisibleByPToEBound2D}}\label{subsubsec3P7P10}

\begin{itemize}
\item Syntax: {\tt NrElementOrdersSharplyDivisibleByPToEBound2D(d,Q,e)}, $d\in\IN^+$, $d\geqslant2$, $Q=p^{2f}$ a prime power square, $e\in\{0,1,\ldots,1+\lceil\log_p(2d-3)\rceil\}$.
\item Output: an upper bound $\overline{\omicron_{\ssrm}^{(e)}}({^2}D_d(Q))$ on the number $\omicron_{\ssrm}^{(e)}({^2}D_d(Q))$ of element orders $o$ in ${^2}D_d(Q)$ such that the largest power of $p$ dividing $o$ is $p^e$.
\item Comments: For an explanation, see \cite[Case (5) of the proof of Theorem 1.1.3(5) in Subsection 3.3]{BGP19a}.
\end{itemize}

\subsubsection{{\tt NrElementOrdersBound2D}}\label{subsubsec3P7P11}

\begin{itemize}
\item Syntax: {\tt NrElementOrdersBound2D(d,Q,l)}, $d\in\IN^+$, $d\geqslant2$, $Q=q^2$ for a prime power $q$, $l\in\IN^+$ ($l$ is interpreted as a \enquote{quality level}).
\item Output: If $l\notin\{1,2,3\}$, it returns the message \enquote{This quality level is not available. Please set the quality level to 1, 2 or 3.} as a string. Otherwise, it returns a certain upper bound $\overline{\omicron}^{(l)}({^2}D_d(Q))$ on the number of element orders $\omicron({^2}D_d(Q))$. The larger $l$ is, the better is that bound, but also the more costly it is to compute.
\item Comments: Write $q=p^f$. For $l=3$, this function returns the sum of the numbers returned by {\tt NrElementOrdersSharplyDivisibleByPToEBound2D(d,Q,e)} where $e=0,1,\ldots,1+\lfloor\log_p(2d-3)\rfloor$. For $l\in\{1,2\}$, the output is $\overline{\omicron_{\ssrm}}^{(l)}({^2}D_d(Q))\cdot(1+\lceil\log_p(2d-2)\rceil)$, the product of an upper bound on $\omicron_{\ssrm}({^2}D_d(Q))$ and the exact number of element orders in ${^2}D_d(Q)$ that are powers of $p$.
\end{itemize}

\subsubsection{{\tt epsilonQBound2D}}\label{subsubsec3P7P12}

\begin{itemize}
\item Syntax: {\tt epsilonQBound2D(d,Q,l1,l2)}, $d\in\IN^+$, $d\geqslant2$, $Q=q^2$ for a prime power $q$, $l_1,l_2\in\IN^+$ (the $l_i$ are interpreted as \enquote{quality levels}).
\item Output: If $(l_1,l_2)\notin\{1,2\}\times\{1,2,3\}$, then it returns the message \enquote{This combination of quality levels is not available. Please set the quality levels to $(1,1)$, $(1,2)$, $(1,3)$, $(2,1)$, $(2,2)$ or $(2,3)$.} as a string. Otherwise, it returns the lower bound
\[
\underline{\epsilon_{\q}}^{(l_1,l_2)}({^2}D_d(Q)):=\frac{\log\log{(\frac{\underline{\omega}^{(l_1)}({^2}D_d(Q))}{\overline{\omicron}^{(l_2)}({^2}D_d(Q)}+3)}}{\log\log{|{^2}D_d(Q)|}}
\]
on $\epsilon_{\q}({^2}D_d(Q))$ as a decimal floating-point number.
\end{itemize}

\subsubsection{{\tt epsilonQBound2Dd4}}\label{subsubsec3P7P13}

\begin{itemize}
\item Syntax: {\tt epsilonQBound2Dd4(d)}, $d\in\IN^+$ (the intended range of use is $d\in\{54,\ldots,90\}$).
\item Output: the real number
\[
\frac{\log\log{(\frac{\underline{\omega}^{(2)}({^2}D_d(4))}{\omicron_{\ssrm}({^2}D_{90}(4))\cdot(1+\lceil\log_2(2d-2)\rceil)}+3)}}{\log\log{|{^2}D_d(4)|}},
\]
as a decimal floating-point number.
\item Comments: If $d\leq90$, the output is a lower bound on $\epsilon_{\q}({^2}D_d(4))$. It is used in \cite[Case (3) of the proof of Theorem 1.1.3(5) in Subsection 3.3]{BGP19a}.
\end{itemize}

\subsection{Functions for groups of type \texorpdfstring{${^2}B_2$}{2B2}}\label{subsec3P8}

\subsubsection{{\tt Order2B2}}\label{subsubsec3P8P1}

\begin{itemize}
\item Syntax: {\tt Order2B2(Q)}, $Q=2^{2k+1}$ for some $k\in\IN^+$.
\item Output: the group order $|{^2}B_2(Q)|$.
\end{itemize}

\subsubsection{{\tt LogLogOrder2B2}}\label{subsubsec3P8P2}

\begin{itemize}
\item Syntax: {\tt LogLogOrder2B2(Q)}, $Q=2^{2k+1}$ for some $k\in\IN^+$.
\item Output: the real number $\log\log{|{^2}B_2(Q)|}$ as a decimal floating-point number.
\item Comments: This does \emph{not} just call {\tt Log(Log(Float(Order2B2(Q))))}, because applying {\tt Float} to too large integers does not result in an appropriate floating-point representation of those integers, but in the output {\tt inf}, which is treated like infinity.
\end{itemize}

\subsubsection{{\tt OrderOut2B2}}\label{subsubsec3P8P3}

\begin{itemize}
\item Syntax: {\tt OrderOut2B2(Q)}, $Q=2^{2k+1}$ for some $k\in\IN^+$.
\item Output: the outer automorphism group order $|\Out({^2}B_2(Q))|$.
\end{itemize}

\subsubsection{{\tt NrConjugacyClasses2B2}}\label{subsec3P8P4}

\begin{itemize}
\item Syntax: {\tt NrConjugacyClasses2B2(Q)}, $Q=2^{2k+1}$ for some $k\in\IN^+$.
\item Output: the conjugacy class number $\k({^2}B_2(Q))$.
\item Comments: See L{\"u}beck's database \cite{Lue}.
\end{itemize}

\subsubsection{{\tt NrAutOrbitsBound2B2}}\label{subsubsec3P8P5}

\begin{itemize}
\item Syntax: {\tt NrAutOrbitsBound2B2(Q)}, $Q=2^{2k+1}$ for some $k\in\IN^+$.
\item Output: the lower bound $\underline{\omega}({^2}B_2(Q)):=\lceil\frac{\k({^2}B_2(Q))}{|\Out({^2}B_2(Q))|}\rceil$ on $\omega({^2}B_2(Q))$.
\end{itemize}

\subsubsection{{\tt epsilonOmegaBound2B2}}\label{subsubsec3P8P6}

\begin{itemize}
\item Syntax: {\tt epsilonOmegaBound2B2(Q)}, $Q=2^{2k+1}$ for some $k\in\IN^+$.
\item Output: the lower bound $\underline{\epsilon_{\omega}}({^2}B_2(Q)):=\frac{\log\log{\underline{\omega}({^2}B_2(Q))}}{\log\log{|{^2}B_2(Q)|}}$ on $\epsilon_{\omega}({^2}B_2(Q))$.
\end{itemize}

\subsubsection{{\tt SemisimpleElementOrders2B2}}\label{subsubsec3P8P7}

\begin{itemize}
\item Syntax: {\tt SemisimpleElementOrders2B2(Q)}, $Q=2^{2k+1}$ for some $k\in\IN^+$.
\item Output: the set of semisimple (viz., odd) element orders in ${^2}B_2(Q)$.
\item Comments: Based on \cite[Theorem 2]{Shi92a}.
\end{itemize}

\subsubsection{{\tt NrSemisimpleElementOrders2B2}}\label{subsubsec3P8P8}

\begin{itemize}
\item Syntax: {\tt NrSemisimpleElementOrders2B2(Q)}, $Q=2^{2k+1}$ for some $k\in\IN^+$.
\item Output: the number $\omicron_{\ssrm}({^2}B_2(Q))$ of semisimple (viz., odd) element orders in ${^2}B_2(Q)$.
\item Comments: Simply calls {\tt Length(SemisimpleElementOrders2B2(Q))}.
\end{itemize}

\subsubsection{{\tt ElementOrders2B2}}\label{subsubsec3P8P9}

\begin{itemize}
\item Syntax: {\tt ElementOrders2B2(Q)}, $Q=2^{2k+1}$ for some $k\in\IN^+$.
\item Output: the set of element orders in ${^2}B_2(Q)$.
\item Comments: Based on \cite[Theorem 2]{Shi92a}.
\end{itemize}

\subsubsection{{\tt NrElementOrders2B2}}\label{subsubsec3P8P10}

\begin{itemize}
\item Syntax: {\tt NrElementOrders2B2(Q)}, $Q=2^{2k+1}$ for some $k\in\IN^+$.
\item Output: the number $\omicron({^2}B_2(Q))$ of element orders in ${^2}B_2(Q)$.
\item Comments: Simply calls {\tt Length(ElementOrders2B2(Q))}.
\end{itemize}

\subsubsection{{\tt epsilonQBound2B2}}\label{subsubsec3P8P11}

\begin{itemize}
\item Syntax: {\tt epsilonQBound2B2(Q)}, $Q=2^{2k+1}$ for some $k\in\IN^+$.
\item Output: the lower bound $\underline{\epsilon_{\q}}({^2}B_2(Q)):=\frac{\log\log{(\underline{\omega}({^2}B_2(Q))/\omicron({^2}B_2(Q))+3)}}{\log\log{|{^2}B_2(Q)|}}$ on $\epsilon_{\q}({^2}B_2(Q))$.
\end{itemize}

\subsection{Functions for groups of type \texorpdfstring{$G_2$}{G2}}\label{subsec3P9}

\subsubsection{{\tt OrderG2}}\label{subsubsec3P9P1}

\begin{itemize}
\item Syntax: {\tt OrderG2(q)}, $q$ a prime power.
\item Output: the group order $|G_2(q)|$.
\end{itemize}

\subsubsection{{\tt LogLogOrderG2}}\label{subsubsec3P9P2}

\begin{itemize}
\item Syntax: {\tt LogLogOrderG2(q)}, $q$ a prime power.
\item Output: the real number $\log\log{|G_2(q)|}$ as a decimal floating-point number.
\item Comments: This does \emph{not} just call {\tt Log(Log(Float(OrderG2(q))))}, because applying {\tt Float} to too large integers does not result in an appropriate floating-point representation of those integers, but in the output {\tt inf}, which is treated like infinity.
\end{itemize}

\subsubsection{{\tt OrderOutG2}}\label{subsubsec3P9P3}

\begin{itemize}
\item Syntax: {\tt OrderOutG2(q)}, $q$ a prime power.
\item Output: the outer automorphism group order $|\Out(G_2(q))|$.
\end{itemize}

\subsubsection{{\tt NrConjugacyClassesG2}}\label{subsec3P9P4}

\begin{itemize}
\item Syntax: {\tt NrConjugacyClassesG2(q)}, $q$ a prime power.
\item Output: the conjugacy class number $\k(G_2(q))$.
\item Comments: See L{\"u}beck's database \cite{Lue}.
\end{itemize}

\subsubsection{{\tt NrAutOrbitsBoundG2}}\label{subsubsec3P9P5}

\begin{itemize}
\item Syntax: {\tt NrAutOrbitsBoundG2(q)}, $q$ a prime power.
\item Output: the lower bound $\underline{\omega}(G_2(q)):=\lceil\frac{\k(G_2(q))}{|\Out(G_2(q))|}\rceil$ on $\omega(G_2(q))$.
\end{itemize}

\subsubsection{{\tt epsilonOmegaBoundG2}}\label{subsubsec3P9P6}

\begin{itemize}
\item Syntax: {\tt epsilonOmegaBoundG2(q)}, $q$ a prime power.
\item Output: the lower bound $\underline{\epsilon_{\omega}}(G_2(q)):=\frac{\log\log{\underline{\omega}(G_2(q))}}{\log\log{|G_2(q)|}}$ on $\epsilon_{\omega}(G_2(q))$.
\end{itemize}

\subsubsection{{\tt CoxeterNoG2}}\label{subsubsec3P9P7}

\begin{itemize}
\item Syntax: {\tt CoxeterNoG2()}.
\item Output: the Coxeter number $h(G_2)=6$.
\end{itemize}

\subsubsection{{\tt SemisimpleElementOrdersG2}}\label{subsubsec3P9P8}

\begin{itemize}
\item Syntax: {\tt SemisimpleElementOrdersG2(q)}, $q=p^f$ a prime power.
\item Output: the set of semisimple (viz., not divisible by $p$) element orders in $G_2(q)$.
\item Comments: Based on \cite[Lemma 1.4]{VS13a}.
\end{itemize}

\subsubsection{{\tt NrSemisimpleElementOrdersG2}}\label{subsubsec3P9P9}

\begin{itemize}
\item Syntax: {\tt NrSemisimpleElementOrdersG2(q)}, $q=p^f$ a prime power.
\item Output: the number $\omicron_{\ssrm}(G_2(q))$ of semisimple (viz., not divisible by $p$) element orders in $G_2(q)$.
\item Comments: Simply calls {\tt Length(SemisimpleElementOrdersG2(q))}.
\end{itemize}

\subsubsection{{\tt ElementOrdersG2}}\label{subsubsec3P9P10}

\begin{itemize}
\item Syntax: {\tt ElementOrdersG2(q)}, $q$ a prime power.
\item Output: the set of element orders in $G_2(q)$.
\item Comments: Based on \cite[Lemma 1.4]{VS13a}.
\end{itemize}

\subsubsection{{\tt NrElementOrdersG2}}\label{subsubsec3P9P11}

\begin{itemize}
\item Syntax: {\tt NrElementOrdersG2(q)}, $q$ a prime power.
\item Output: the number $\omicron(G_2(q))$ of element orders in $G_2(q)$.
\item Comments: Simply calls {\tt Length(ElementOrdersG2(q))}.
\end{itemize}

\subsubsection{{\tt epsilonQBoundG2}}\label{subsubsec3P9P12}

\begin{itemize}
\item Syntax: {\tt epsilonQBoundG2(q)}, $q$ a prime power.
\item Output: the lower bound $\underline{\epsilon_{\q}}(G_2(q)):=\frac{\log\log{(\underline{\omega}(G_2(q))/\omicron(G_2(q))+3)}}{\log\log{|G_2(q)|}}$ on $\epsilon_{\q}(G_2(q))$.
\end{itemize}

\subsection{Functions for groups of type \texorpdfstring{${^2}G_2$}{2G2}}\label{subsec3P10}

\subsubsection{{\tt Order2G2}}\label{subsubsec3P10P1}

\begin{itemize}
\item Syntax: {\tt Order2G2(Q)}, $Q=3^{2k+1}$ for some $k\in\IN^+$.
\item Output: the group order $|{^2}G_2(Q)|$.
\end{itemize}

\subsubsection{{\tt LogLogOrder2G2}}\label{subsubsec3P10P2}

\begin{itemize}
\item Syntax: {\tt LogLogOrder2G2(Q)}, $Q=3^{2k+1}$ for some $k\in\IN^+$.
\item Output: the real number $\log\log{|{^2}G_2(Q)|}$ as a decimal floating-point number.
\item Comments: This does \emph{not} just call {\tt Log(Log(Float(Order2G2(Q))))}, because applying {\tt Float} to too large integers does not result in an appropriate floating-point representation of those integers, but in the output {\tt inf}, which is treated like infinity.
\end{itemize}

\subsubsection{{\tt OrderOut2G2}}\label{subsubsec3P10P3}

\begin{itemize}
\item Syntax: {\tt OrderOut2G2(Q)}, $Q=3^{2k+1}$ for some $k\in\IN^+$.
\item Output: the outer automorphism group order $|\Out({^2}G_2(Q))|$.
\end{itemize}

\subsubsection{{\tt NrConjugacyClasses2G2}}\label{subsec3P10P4}

\begin{itemize}
\item Syntax: {\tt NrConjugacyClasses2G2(Q)}, $Q=3^{2k+1}$ for some $k\in\IN^+$.
\item Output: the conjugacy class number $\k({^2}G_2(Q))$.
\item Comments: See L{\"u}beck's database \cite{Lue}.
\end{itemize}

\subsubsection{{\tt NrAutOrbitsBound2G2}}\label{subsubsec3P10P5}

\begin{itemize}
\item Syntax: {\tt NrAutOrbitsBound2G2(Q)}, $Q=3^{2k+1}$ for some $k\in\IN^+$.
\item Output: the lower bound $\underline{\omega}({^2}G_2(Q)):=\lceil\frac{\k({^2}G_2(Q))}{|\Out({^2}G_2(Q))|}\rceil$ on $\omega({^2}G_2(Q))$.
\end{itemize}

\subsubsection{{\tt epsilonOmegaBound2G2}}\label{subsubsec3P10P6}

\begin{itemize}
\item Syntax: {\tt epsilonOmegaBound2G2(Q)}, $Q=3^{2k+1}$ for some $k\in\IN^+$.
\item Output: the lower bound $\underline{\epsilon_{\omega}}({^2}G_2(Q)):=\frac{\log\log{\underline{\omega}({^2}G_2(Q))}}{\log\log{|{^2}G_2(Q)|}}$ on $\epsilon_{\omega}({^2}G_2(Q))$.
\end{itemize}

\subsubsection{{\tt SemisimpleElementOrders2G2}}\label{subsubsec3P10P7}

\begin{itemize}
\item Syntax: {\tt SemisimpleElementOrders2G2(Q)}, $Q=3^{2k+1}$ for some $k\in\IN^+$.
\item Output: the set of semisimple (viz., not divisible by $3$) element orders in ${^2}G_2(Q)$.
\item Comments: Based on \cite[Lemma 4]{BS93a}.
\end{itemize}

\subsubsection{{\tt NrSemisimpleElementOrders2G2}}\label{subsubsec3P10P8}

\begin{itemize}
\item Syntax: {\tt NrSemisimpleElementOrders2G2(Q)}, $Q=3^{2k+1}$ for some $k\in\IN^+$.
\item Output: the number $\omicron_{\ssrm}({^2}G_2(Q))$ of semisimple (viz., not divisible by $3$) element orders in ${^2}G_2(Q)$.
\item Comments: Simply calls {\tt Length(SemisimpleElementOrders2G2(Q))}.
\end{itemize}

\subsubsection{{\tt ElementOrders2G2}}\label{subsubsec3P10P9}

\begin{itemize}
\item Syntax: {\tt ElementOrders2G2(Q)}, $Q=3^{2k+1}$ for some $k\in\IN^+$.
\item Output: the set of element orders in ${^2}G_2(Q)$.
\item Comments: Based on \cite[Lemma 4]{BS93a}.
\end{itemize}

\subsubsection{{\tt NrElementOrders2G2}}\label{subsubsec3P10P10}

\begin{itemize}
\item Syntax: {\tt NrElementOrders2G2(Q)}, $Q=3^{2k+1}$ for some $k\in\IN^+$.
\item Output: the number $\omicron({^2}G_2(Q))$ of element orders in ${^2}G_2(Q)$.
\item Comments: Simply calls {\tt Length(ElementOrders2G2(Q))}.
\end{itemize}

\subsubsection{{\tt epsilonQBound2G2}}\label{subsubsec3P10P11}

\begin{itemize}
\item Syntax: {\tt epsilonQBound2G2(Q)}, $Q=3^{2k+1}$ for some $k\in\IN^+$.
\item Output: the lower bound $\underline{\epsilon_{\q}}({^2}G_2(Q)):=\frac{\log\log{(\underline{\omega}({^2}G_2(Q))/\omicron({^2}G_2(Q))+3)}}{\log\log{|{^2}G_2(Q)|}}$ on $\epsilon_{\q}({^2}G_2(Q))$.
\end{itemize}

\subsection{Functions for groups of type \texorpdfstring{${^3}D_4$}{3D4}}\label{subsec3P11}

\subsubsection{{\tt Order3D4}}\label{subsubsec3P11P1}

\begin{itemize}
\item Syntax: {\tt Order3D4(Q)}, $Q=q^3$ for some prime power $q$.
\item Output: the group order $|{^3}D_4(q)|$.
\end{itemize}

\subsubsection{{\tt LogLogOrder3D4}}\label{subsubsec3P11P2}

\begin{itemize}
\item Syntax: {\tt LogLogOrder3D4(Q)}, $Q=q^3$ for some prime power $q$.
\item Output: the real number $\log\log{|{^3}D_4(Q)|}$ as a decimal floating-point number.
\item Comments: This does \emph{not} just call {\tt Log(Log(Float(Order3D4(Q))))}, because applying {\tt Float} to too large integers does not result in an appropriate floating-point representation of those integers, but in the output {\tt inf}, which is treated like infinity.
\end{itemize}

\subsubsection{{\tt OrderOut3D4}}\label{subsubsec3P11P3}

\begin{itemize}
\item Syntax: {\tt OrderOut3D4(Q)}, $Q=q^3$ for some prime power $q$.
\item Output: the outer automorphism group order $|\Out({^3}D_4(q))|$.
\end{itemize}

\subsubsection{{\tt NrConjugacyClasses3D4}}\label{subsec3P11P4}

\begin{itemize}
\item Syntax: {\tt NrConjugacyClasses3D4(Q)}, $Q=q^3$ for some prime power $q$.
\item Output: the conjugacy class number $\k({^3}D_4(Q))$.
\item Comments: See L{\"u}beck's database \cite{Lue}.
\end{itemize}

\subsubsection{{\tt NrAutOrbitsBound3D4}}\label{subsubsec3P11P5}

\begin{itemize}
\item Syntax: {\tt NrAutOrbitsBound3D4(Q)}, $Q=q^3$ for some prime power $q$.
\item Output: the lower bound $\underline{\omega}({^3}D_4(Q)):=\lceil\frac{\k({^3}D_4(Q))}{|\Out({^3}D_4(Q))|}\rceil$ on $\omega({^3}D_4(Q))$.
\end{itemize}

\subsubsection{{\tt epsilonOmegaBound3D4}}\label{subsubsec3P11P6}

\begin{itemize}
\item Syntax: {\tt epsilonOmegaBound3D4(Q)}, $Q=q^3$ for some prime power $q$.
\item Output: the lower bound $\underline{\epsilon_{\omega}}({^3}D_4(Q)):=\frac{\log\log{\underline{\omega}({^3}D_4(Q))}}{\log\log{|{^3}D_4(Q)|}}$ on $\epsilon_{\omega}({^3}D_4(Q))$.
\end{itemize}

\subsubsection{{\tt SemisimpleElementOrders3D4}}\label{subsubsec3P11P7}

\begin{itemize}
\item Syntax: {\tt SemisimpleElementOrders3D4(Q)}, $Q=p^{3f}$ some prime power cube.
\item Output: the set of semisimple (viz., not divisible by $3$) element orders in ${^3}D_4(Q)$.
\item Comments: Based on \cite[Theorem 3.2]{GZ16a}.
\end{itemize}

\subsubsection{{\tt NrSemisimpleElementOrders3D4}}\label{subsubsec3P11P8}

\begin{itemize}
\item Syntax: {\tt NrSemisimpleElementOrders3D4(Q)}, $Q=p^{3f}$ some prime power cube.
\item Output: the number $\omicron_{\ssrm}({^3}D_4(Q))$ of semisimple (viz., not divisible by $p$) element orders in ${^3}D_4(Q)$.
\item Comments: Simply calls {\tt Length(SemisimpleElementOrders3D4(Q))}.
\end{itemize}

\subsubsection{{\tt ElementOrders3D4}}\label{subsubsec3P11P9}

\begin{itemize}
\item Syntax: {\tt ElementOrders3D4(Q)}, $Q=q^3$ for some prime power cube.
\item Output: the set of element orders in ${^3}D_4(Q)$.
\item Comments: Based on \cite[Theorem 3.2]{GZ16a}.
\end{itemize}

\subsubsection{{\tt NrElementOrders3D4}}\label{subsubsec3P11P10}

\begin{itemize}
\item Syntax: {\tt NrElementOrders3D4(Q)}, $Q=q^3$ for some prime power $q$.
\item Output: the number $\omicron({^3}D_4(Q))$ of element orders in ${^3}D_4(Q)$.
\item Comments: Simply calls {\tt Length(ElementOrders3D4(Q))}.
\end{itemize}

\subsubsection{{\tt epsilonQBound3D4}}\label{subsubsec3P11P11}

\begin{itemize}
\item Syntax: {\tt epsilonQBound3D4(Q)}, $Q=q^3$ for some prime power $q$.
\item Output: the lower bound $\underline{\epsilon_{\q}}({^3}D_4(Q)):=\frac{\log\log{(\underline{\omega}({^3}D_4(Q))/\omicron({^3}D_4(Q))+3)}}{\log\log{|{^3}D_4(Q)|}}$ on $\epsilon_{\q}({^3}D_4(Q))$.
\end{itemize}

\subsection{Functions for groups of type \texorpdfstring{$F_4$}{F4}}\label{subsec3P12}

\subsubsection{{\tt OrderF4}}\label{subsubsec3P12P1}

\begin{itemize}
\item Syntax: {\tt OrderF4(q)}, $q$ a prime power.
\item Output: the group order $|F_4(q)|$.
\end{itemize}

\subsubsection{{\tt LogLogOrderF4}}\label{subsubsec3P12P2}

\begin{itemize}
\item Syntax: {\tt LogLogOrderF4(q)}, $q$ a prime power.
\item Output: the real number $\log\log{|F_4(q)|}$ as a decimal floating-point number.
\item Comments: This does \emph{not} just call {\tt Log(Log(Float(OrderF4(q))))}, because applying {\tt Float} to too large integers does not result in an appropriate floating-point representation of those integers, but in the output {\tt inf}, which is treated like infinity.
\end{itemize}

\subsubsection{{\tt OrderOutF4}}\label{subsubsec3P12P3}

\begin{itemize}
\item Syntax: {\tt OrderOutF4(q)}, $q$ a prime power.
\item Output: the outer automorphism group order $|\Out(F_4(q))|$.
\end{itemize}

\subsubsection{{\tt NrConjugacyClassesF4}}\label{subsec3P12P4}

\begin{itemize}
\item Syntax: {\tt NrConjugacyClassesF4(q)}, $q$ a prime power.
\item Output: the conjugacy class number $\k(F_4(q))$.
\item Comments: See L{\"u}beck's database \cite{Lue}.
\end{itemize}

\subsubsection{{\tt NrAutOrbitsBoundF4}}\label{subsubsec3P12P5}

\begin{itemize}
\item Syntax: {\tt NrAutOrbitsBoundF4(q)}, $q$ a prime power.
\item Output: the lower bound $\underline{\omega}(F_4(q)):=\lceil\frac{\k(F_4(q))}{|\Out(F_4(q))|}\rceil$ on $\omega(F_4(q))$.
\end{itemize}

\subsubsection{{\tt epsilonOmegaBoundF4}}\label{subsubsec3P12P6}

\begin{itemize}
\item Syntax: {\tt epsilonOmegaBoundF4(q)}, $q$ a prime power.
\item Output: the lower bound $\underline{\epsilon_{\omega}}(F_4(q)):=\frac{\log\log{\underline{\omega}(F_4(q))}}{\log\log{|F_4(q)|}}$ on $\epsilon_{\omega}(F_4(q))$.
\end{itemize}

\subsubsection{{\tt CoxeterNoF4}}\label{subsubsec3P12P7}

\begin{itemize}
\item Syntax: {\tt CoxeterNoF4()}.
\item Output: the Coxeter number $h(F_4)=12$.
\end{itemize}

\subsubsection{{\tt SemisimpleElementOrdersF4}}\label{subsubsec3P12P8}

\begin{itemize}
\item Syntax: {\tt SemisimpleElementOrdersF4(q)}, $q=p^f$ a prime power.
\item Output: the set of semisimple (viz., not divisible by $p$) element orders in $F_4(q)$.
\item Comments: Based on \cite[Theorem 3.1]{GZ16a}.
\end{itemize}

\subsubsection{{\tt NrSemisimpleElementOrdersF4}}\label{subsubsec3P12P9}

\begin{itemize}
\item Syntax: {\tt NrSemisimpleElementOrdersF4(q)}, $q=p^f$ a prime power.
\item Output: the number $\omicron_{\ssrm}(F_4(q))$ of semisimple (viz., not divisible by $p$) element orders in $F_4(q)$.
\item Comments: Simply calls {\tt Length(SemisimpleElementOrdersF4(q))}.
\end{itemize}

\subsubsection{{\tt ElementOrdersF4}}\label{subsubsec3P12P10}

\begin{itemize}
\item Syntax: {\tt ElementOrdersF4(q)}, $q$ a prime power.
\item Output: the set of element orders in $F_4(q)$.
\item Comments: Based on \cite[Theorem 3.1]{GZ16a}.
\end{itemize}

\subsubsection{{\tt NrElementOrdersF4}}\label{subsubsec3P12P11}

\begin{itemize}
\item Syntax: {\tt NrElementOrdersF4(q)}, $q$ a prime power.
\item Output: the number $\omicron(F_4(q))$ of element orders in $F_4(q)$.
\item Comments: Simply calls {\tt Length(ElementOrdersF4(q))}.
\end{itemize}

\subsubsection{{\tt epsilonQBoundF4}}\label{subsubsec3P12P12}

\begin{itemize}
\item Syntax: {\tt epsilonQBoundF4(q)}, $q$ a prime power.
\item Output: the lower bound $\underline{\epsilon_{\q}}(F_4(q)):=\frac{\log\log{(\underline{\omega}(F_4(q))/\omicron(F_4(q))+3)}}{\log\log{|F_4(q)|}}$ on $\epsilon_{\q}(F_4(q))$.
\end{itemize}

\subsection{Functions for groups of type \texorpdfstring{${^2}F_4$}{2F4}}\label{subsec3P13}

\subsubsection{{\tt Order2F4}}\label{subsubsec3P13P1}

\begin{itemize}
\item Syntax: {\tt Order2F4(Q)}, $Q=2^{2k+1}$ for some $k\in\IN^+$.
\item Output: the group order $|{^2}F_4(Q)|$.
\end{itemize}

\subsubsection{{\tt LogLogOrder2F4}}\label{subsubsec3P13P2}

\begin{itemize}
\item Syntax: {\tt LogLogOrder2F4(Q)}, $Q=2^{2k+1}$ for some $k\in\IN^+$.
\item Output: the real number $\log\log{|{^2}F_4(Q)|}$ as a decimal floating-point number.
\item Comments: This does \emph{not} just call {\tt Log(Log(Float(Order2F4(Q))))}, because applying {\tt Float} to too large integers does not result in an appropriate floating-point representation of those integers, but in the output {\tt inf}, which is treated like infinity.
\end{itemize}

\subsubsection{{\tt OrderOut2F4}}\label{subsubsec3P13P3}

\begin{itemize}
\item Syntax: {\tt OrderOut2F4(Q)}, $Q=2^{2k+1}$ for some $k\in\IN^+$.
\item Output: the outer automorphism group order $|\Out({^2}F_4(Q))|$.
\end{itemize}

\subsubsection{{\tt NrConjugacyClasses2F4}}\label{subsec3P13P4}

\begin{itemize}
\item Syntax: {\tt NrConjugacyClasses2F4(Q)}, $Q=2^{2k+1}$ for some $k\in\IN^+$.
\item Output: the conjugacy class number $\k({^2}F_4(Q))$.
\item Comments: See L{\"u}beck's database \cite{Lue}.
\end{itemize}

\subsubsection{{\tt NrAutOrbitsBound2F4}}\label{subsubsec3P13P5}

\begin{itemize}
\item Syntax: {\tt NrAutOrbitsBound2F4(Q)}, $Q=2^{2k+1}$ for some $k\in\IN^+$.
\item Output: the lower bound $\underline{\omega}({^2}F_4(Q)):=\lceil\frac{\k({^2}F_4(Q))}{|\Out({^2}F_4(Q))|}\rceil$ on $\omega({^2}F_4(Q))$.
\end{itemize}

\subsubsection{{\tt epsilonOmegaBound2F4}}\label{subsubsec3P13P6}

\begin{itemize}
\item Syntax: {\tt epsilonOmegaBound2F4(Q)}, $Q=2^{2k+1}$ for some $k\in\IN^+$.
\item Output: the lower bound $\underline{\epsilon_{\omega}}({^2}F_4(Q)):=\frac{\log\log{\underline{\omega}({^2}F_4(Q))}}{\log\log{|{^2}F_4(Q)|}}$ on $\epsilon_{\omega}({^2}F_4(Q))$.
\end{itemize}

\subsubsection{{\tt SemisimpleElementOrders2F4}}\label{subsubsec3P13P7}

\begin{itemize}
\item Syntax: {\tt SemisimpleElementOrders2F4(Q)}, $Q=2^{2k+1}$ for some $k\in\IN^+$.
\item Output: the set of semisimple (viz., odd) element orders in ${^2}F_4(Q)$.
\item Comments: Based on \cite[Lemma 3]{DS99a}.
\end{itemize}

\subsubsection{{\tt NrSemisimpleElementOrders2F4}}\label{subsubsec3P13P8}

\begin{itemize}
\item Syntax: {\tt NrSemisimpleElementOrders2F4(Q)}, $Q=2^{2k+1}$ for some $k\in\IN^+$.
\item Output: the number $\omicron_{\ssrm}({^2}F_4(Q))$ of semisimple (viz., odd) element orders in ${^2}F_4(Q)$.
\item Comments: Simply calls {\tt Length(SemisimpleElementOrders2F4(Q))}.
\end{itemize}

\subsubsection{{\tt ElementOrders2F4}}\label{subsubsec3P13P9}

\begin{itemize}
\item Syntax: {\tt ElementOrders2F4(Q)}, $Q=2^{2k+1}$ for some $k\in\IN^+$.
\item Output: the set of element orders in ${^2}F_4(Q)$.
\item Comments: Based on \cite[Lemma 3]{DS99a}.
\end{itemize}

\subsubsection{{\tt NrElementOrders2F4}}\label{subsubsec3P13P10}

\begin{itemize}
\item Syntax: {\tt NrElementOrders2F4(Q)}, $Q=2^{2k+1}$ for some $k\in\IN^+$.
\item Output: the number $\omicron({^2}F_4(Q))$ of element orders in ${^2}F_4(Q)$.
\item Comments: Simply calls {\tt Length(ElementOrdersF4(q))}.
\end{itemize}

\subsubsection{{\tt epsilonQBound2F4}}\label{subsubsec3P13P11}

\begin{itemize}
\item Syntax: {\tt epsilonQBound2F4(Q)}, $Q=2^{2k+1}$ for some $k\in\IN^+$.
\item Output: the lower bound $\underline{\epsilon_{\q}}({^2}F_4(Q)):=\frac{\log\log{(\underline{\omega}({^2}F_4(Q))/\omicron({^2}F_4(Q))+3)}}{\log\log{|{^2}F_4(Q)|}}$ on $\epsilon_{\q}({^2}F_4(Q))$.
\end{itemize}

\subsection{Functions for groups of type \texorpdfstring{$E_6$}{E6}}\label{subsec3P14}

\subsubsection{{\tt OrderE6}}\label{subsubsec3P14P1}

\begin{itemize}
\item Syntax: {\tt OrderE6(q)}, $q$ a prime power.
\item Output: the group order $|E_6(q)|$.
\end{itemize}

\subsubsection{{\tt LogLogOrderE6}}\label{subsubsec3P14P2}

\begin{itemize}
\item Syntax: {\tt LogLogOrderE6(q)}, $q$ a prime power.
\item Output: the real number $\log\log{|E_6(q)|}$ as a decimal floating-point number.
\item Comments: This does \emph{not} just call {\tt Log(Log(Float(OrderE6(q))))}, because applying {\tt Float} to too large integers does not result in an appropriate floating-point representation of those integers, but in the output {\tt inf}, which is treated like infinity.
\end{itemize}

\subsubsection{{\tt OrderOutE6}}\label{subsubsec3P14P3}

\begin{itemize}
\item Syntax: {\tt OrderOutE6(q)}, $q$ a prime power.
\item Output: the outer automorphism group order $|\Out(E_6(q))|$.
\end{itemize}

\subsubsection{{\tt NrConjugacyClassesBoundE6}}\label{subsec3P14P4}

\begin{itemize}
\item Syntax: {\tt NrConjugacyClassesBoundE6(q)}, $q$ a prime power.
\item Output: the lower bound $\underline{\k}(E_6(q)):=\lceil\frac{\k(\Inndiag(E_6(q)))}{|\Outdiag(E_6(q))|}\rceil$ on $\k(E_6(q))$.
\item Comments: For $\k(\Inndiag(E_6(q)))$ see L{\"u}beck's database \cite{Lue}.
\end{itemize}

\subsubsection{{\tt NrAutOrbitsBoundE6}}\label{subsubsec3P14P5}

\begin{itemize}
\item Syntax: {\tt NrAutOrbitsBoundE6(q)}, $q$ a prime power.
\item Output: the lower bound $\underline{\omega}(E_6(q)):=\lceil\frac{\underline{\k}(E_6(q))}{|\Out(E_6(q))|}\rceil$ on $\omega(E_6(q))$.
\end{itemize}

\subsubsection{{\tt epsilonOmegaBoundE6}}\label{subsubsec3P14P6}

\begin{itemize}
\item Syntax: {\tt epsilonOmegaBoundE6(q)}, $q$ a prime power.
\item Output: the lower bound $\underline{\epsilon_{\omega}}(E_6(q)):=\frac{\log\log{\underline{\omega}(E_6(q))}}{\log\log{|E_6(q)|}}$ on $\epsilon_{\omega}(E_6(q))$.
\end{itemize}

\subsubsection{{\tt CoxeterNoE6}}\label{subsubsec3P14P7}

\begin{itemize}
\item Syntax: {\tt CoxeterNoE6()}.
\item Output: the Coxeter number $h(E_6)=12$.
\end{itemize}

\subsubsection{{\tt SemisimpleElementOrdersE6}}\label{subsubsec3P14P8}

\begin{itemize}
\item Syntax: {\tt SemisimpleElementOrdersE6(q)}, $q=p^f$ a prime power.
\item Output: the set of semisimple (viz., not divisible by $p$) element orders in $E_6(q)$.
\item Comments: Based on \cite[Theorem 1]{But13a}.
\end{itemize}

\subsubsection{{\tt NrSemisimpleElementOrdersE6}}\label{subsubsec3P14P9}

\begin{itemize}
\item Syntax: {\tt NrSemisimpleElementOrdersE6(q)}, $q=p^f$ a prime power.
\item Output: the number $\omicron_{\ssrm}(E_6(q))$ of semisimple (viz., not divisible by $p$) element orders in $E_6(q)$.
\item Comments: Simply calls {\tt Length(SemisimpleElementOrdersE6(q))}.
\end{itemize}

\subsubsection{{\tt ElementOrdersE6}}\label{subsubsec3P14P10}

\begin{itemize}
\item Syntax: {\tt ElementOrdersE6(q)}, $q$ a prime power.
\item Output: the set of element orders in $E_6(q)$.
\item Comments: Based on \cite[Theorem 1]{But13a}.
\end{itemize}

\subsubsection{{\tt NrElementOrdersE6}}\label{subsubsec3P14P11}

\begin{itemize}
\item Syntax: {\tt NrElementOrdersE6(q)}, $q$ a prime power.
\item Output: the number $\omicron(E_6(q))$ of element orders in $E_6(q)$.
\item Comments: Simply calls {\tt Length(ElementOrdersE6(q))}.
\end{itemize}

\subsubsection{{\tt epsilonQBoundE6}}\label{subsubsec3P14P12}

\begin{itemize}
\item Syntax: {\tt epsilonQBoundE6(q)}, $q$ a prime power.
\item Output: the lower bound $\underline{\epsilon_{\q}}(E_6(q)):=\frac{\log\log{(\underline{\omega}(E_6(q))/\omicron(E_6(q))+3)}}{\log\log{|E_6(q)|}}$ on $\epsilon_{\q}(E_6(q))$.
\end{itemize}

\subsection{Functions for groups of type \texorpdfstring{${^2}E_6$}{2E6}}\label{subsec3P15}

\subsubsection{{\tt Order2E6}}\label{subsubsec3P15P1}

\begin{itemize}
\item Syntax: {\tt Order2E6(Q)}, $Q=q^2$ for a prime power $q$.
\item Output: the group order $|{^2}E_6(q^2)|$.
\end{itemize}

\subsubsection{{\tt LogLogOrder2E6}}\label{subsubsec3P15P2}

\begin{itemize}
\item Syntax: {\tt LogLogOrder2E6(Q)}, $Q=q^2$ for a prime power $q$.
\item Output: the real number $\log\log{|{^2}E_6(Q)|}$ as a decimal floating-point number.
\item Comments: This does \emph{not} just call {\tt Log(Log(Float(Order2E6(Q))))}, because applying {\tt Float} to too large integers does not result in an appropriate floating-point representation of those integers, but in the output {\tt inf}, which is treated like infinity.
\end{itemize}

\subsubsection{{\tt OrderOut2E6}}\label{subsubsec3P15P3}

\begin{itemize}
\item Syntax: {\tt OrderOut2E6(Q)}, $Q=q^2$ for a prime power $q$.
\item Output: the outer automorphism group order $|\Out({^2}E_6(Q))|$.
\end{itemize}

\subsubsection{{\tt NrConjugacyClassesBound2E6}}\label{subsec3P15P4}

\begin{itemize}
\item Syntax: {\tt NrConjugacyClassesBound2E6(Q)}, $Q=q^2$ for a prime power $q$.
\item Output: the lower bound $\underline{\k}({^2}E_6(Q)):=\lceil\frac{\k(\Inndiag({^2}E_6(Q)))}{|\Outdiag({^2}E_6(Q))|}\rceil$ on $\k({^2}E_6(Q))$.
\item Comments: For $\k(\Inndiag({^2}E_6(Q)))$ see L{\"u}beck's database \cite{Lue}.
\end{itemize}

\subsubsection{{\tt NrAutOrbitsBound2E6}}\label{subsubsec3P15P5}

\begin{itemize}
\item Syntax: {\tt NrAutOrbitsBound2E6(Q)}, $Q=q^2$ for a prime power $q$.
\item Output: the lower bound $\underline{\omega}({^2}E_6(Q)):=\lceil\frac{\underline{\k}({^2}E_6(Q))}{|\Out({^2}E_6(Q))|}\rceil$ on $\omega({^2}E_6(Q))$.
\end{itemize}

\subsubsection{{\tt epsilonOmegaBound2E6}}\label{subsubsec3P15P6}

\begin{itemize}
\item Syntax: {\tt epsilonOmegaBound2E6(Q)}, $Q=q^2$ for a prime power $q$.
\item Output: the lower bound $\underline{\epsilon_{\omega}}({^2}E_6(Q)):=\frac{\log\log{\underline{\omega}({^2}E_6(Q))}}{\log\log{|{^2}E_6(Q)|}}$ on $\epsilon_{\omega}({^2}E_6(Q))$.
\end{itemize}

\subsubsection{{\tt SemisimpleElementOrders2E6}}\label{subsubsec3P15P7}

\begin{itemize}
\item Syntax: {\tt SemisimpleElementOrders2E6(Q)}, $Q=p^{2f}$ a prime power square.
\item Output: the set of semisimple (viz., not divisible by $p$) element orders in ${^2}E_6(Q)$.
\item Comments: Based on \cite[Theorem 1]{But13a}.
\end{itemize}

\subsubsection{{\tt NrSemisimpleElementOrders2E6}}\label{subsubsec3P15P8}

\begin{itemize}
\item Syntax: {\tt NrSemisimpleElementOrders2E6(Q)}, $Q=p^{2f}$ a prime power square.
\item Output: the number $\omicron_{\ssrm}({^2}E_6(Q))$ of semisimple (viz., not divisible by $p$) element orders in ${^2}E_6(Q)$.
\item Comments: Simply calls {\tt Length(SemisimpleElementOrders2E6(Q))}.
\end{itemize}

\subsubsection{{\tt ElementOrders2E6}}\label{subsubsec3P15P9}

\begin{itemize}
\item Syntax: {\tt ElementOrders2E6(Q)}, $Q=q^2$ for a prime power $q$.
\item Output: the set of element orders in ${^2}E_6(Q)$.
\item Comments: Based on \cite[Theorem 1]{But13a}.
\end{itemize}

\subsubsection{{\tt NrElementOrders2E6}}\label{subsubsec3P15P10}

\begin{itemize}
\item Syntax: {\tt NrElementOrders2E6(Q)}, $Q=q^2$ for a prime power $q$.
\item Output: the number $\omicron({^2}E_6(Q))$ of element orders in ${^2}E_6(Q)$.
\item Comments: Simply calls {\tt Length(ElementOrders2E6(Q))}.
\end{itemize}

\subsubsection{{\tt epsilonQBound2E6}}\label{subsubsec3P15P11}

\begin{itemize}
\item Syntax: {\tt epsilonQBound2E6(Q)}, $Q=q^2$ for a prime power $q$.
\item Output: the lower bound $\underline{\epsilon_{\q}}({^2}E_6(Q)):=\frac{\log\log{(\underline{\omega}({^2}E_6(Q))/\omicron({^2}E_6(Q))+3)}}{\log\log{|{^2}E_6(Q)|}}$ on $\epsilon_{\q}({^2}E_6(Q))$.
\end{itemize}

\subsection{Functions for groups of type \texorpdfstring{$E_7$}{E7}}\label{subsec3P16}

\subsubsection{{\tt OrderE7}}\label{subsubsec3P16P1}

\begin{itemize}
\item Syntax: {\tt OrderE7(q)}, $q$ a prime power.
\item Output: the group order $|E_7(q)|$.
\end{itemize}

\subsubsection{{\tt LogLogOrderE7}}\label{subsubsec3P16P2}

\begin{itemize}
\item Syntax: {\tt LogLogOrderE7(q)}, $q$ a prime power.
\item Output: the real number $\log\log{|E_7(q)|}$ as a decimal floating-point number.
\item Comments: This does \emph{not} just call {\tt Log(Log(Float(OrderE7(q))))}, because applying {\tt Float} to too large integers does not result in an appropriate floating-point representation of those integers, but in the output {\tt inf}, which is treated like infinity.
\end{itemize}

\subsubsection{{\tt OrderOutE7}}\label{subsubsec3P16P3}

\begin{itemize}
\item Syntax: {\tt OrderOutE7(q)}, $q$ a prime power.
\item Output: the outer automorphism group order $|\Out(E_7(q))|$.
\end{itemize}

\subsubsection{{\tt NrConjugacyClassesBoundE7}}\label{subsec3P16P4}

\begin{itemize}
\item Syntax: {\tt NrConjugacyClassesBoundE7(q)}, $q$ a prime power.
\item Output: the lower bound $\underline{\k}(E_7(q)):=\lceil\frac{\k(\Inndiag(E_7(q)))}{|\Outdiag(E_7(q))|}\rceil$ on $\k(E_7(q))$.
\item Comments: For $\k(\Inndiag(E_7(q)))$ see L{\"u}beck's database \cite{Lue}.
\end{itemize}

\subsubsection{{\tt NrAutOrbitsBoundE7}}\label{subsubsec3P16P5}

\begin{itemize}
\item Syntax: {\tt NrAutOrbitsBoundE7(q)}, $q$ a prime power.
\item Output: the lower bound $\underline{\omega}(E_7(q)):=\lceil\frac{\underline{\k}(E_7(q))}{|\Out(E_7(q))|}\rceil$ on $\omega(E_7(q))$.
\end{itemize}

\subsubsection{{\tt epsilonOmegaBoundE7}}\label{subsubsec3P16P6}

\begin{itemize}
\item Syntax: {\tt epsilonOmegaBoundE7(q)}, $q$ a prime power.
\item Output: the lower bound $\underline{\epsilon_{\omega}}(E_7(q)):=\frac{\log\log{\underline{\omega}(E_7(q))}}{\log\log{|E_7(q)|}}$ on $\epsilon_{\omega}(E_7(q))$.
\end{itemize}

\subsubsection{{\tt CoxeterNoE7}}\label{subsubsec3P16P7}

\begin{itemize}
\item Syntax: {\tt CoxeterNoE7()}.
\item Output: the Coxeter number $h(E_7)=18$.
\end{itemize}

\subsubsection{{\tt SemisimpleElementOrdersE7}}\label{subsubsec3P16P8}

\begin{itemize}
\item Syntax: {\tt SemisimpleElementOrdersE7(q)}, $q=p^f$ a prime power.
\item Output: the set of semisimple (viz., not divisible by $p$) element orders in $E_7(q)$.
\item Comments: Based on \cite[Theorem 2]{But16a}.
\end{itemize}

\subsubsection{{\tt NrSemisimpleElementOrdersE7}}\label{subsubsec3P16P9}

\begin{itemize}
\item Syntax: {\tt NrSemisimpleElementOrdersE7(q)}, $q=p^f$ a prime power.
\item Output: the number $\omicron_{\ssrm}(E_7(q))$ of semisimple (viz., not divisible by $p$) element orders in $E_7(q)$.
\item Comments: Simply calls {\tt Length(SemisimpleElementOrdersE7(q))}.
\end{itemize}

\subsubsection{{\tt ElementOrdersE7}}\label{subsubsec3P16P10}

\begin{itemize}
\item Syntax: {\tt ElementOrdersE7(q)}, $q$ a prime power.
\item Output: the set of element orders in $E_7(q)$.
\item Comments: Based on \cite[Theorem 2]{But16a}.
\end{itemize}

\subsubsection{{\tt NrElementOrdersE7}}\label{subsubsec3P16P11}

\begin{itemize}
\item Syntax: {\tt NrElementOrdersE7(q)}, $q$ a prime power.
\item Output: the number $\omicron(E_7(q))$ of element orders in $E_7(q)$.
\item Comments: Simply calls {\tt Length(ElementOrdersE7(q))}.
\end{itemize}

\subsubsection{{\tt epsilonQBoundE7}}\label{subsubsec3P16P12}

\begin{itemize}
\item Syntax: {\tt epsilonQBoundE7(q)}, $q$ a prime power.
\item Output: the lower bound $\underline{\epsilon_{\q}}(E_7(q)):=\frac{\log\log{(\underline{\omega}(E_7(q))/\omicron(E_7(q))+3)}}{\log\log{|E_7(q)|}}$ on $\epsilon_{\q}(E_7(q))$.
\end{itemize}

\subsection{Functions for groups of type \texorpdfstring{$E_8$}{E8}}\label{subsec3P17}

\subsubsection{{\tt OrderE8}}\label{subsubsec3P17P1}

\begin{itemize}
\item Syntax: {\tt OrderE8(q)}, $q$ a prime power.
\item Output: the group order $|E_8(q)|$.
\end{itemize}

\subsubsection{{\tt LogLogOrderE8}}\label{subsubsec3P17P2}

\begin{itemize}
\item Syntax: {\tt LogLogOrderE8(q)}, $q$ a prime power.
\item Output: the real number $\log\log{|E_8(q)|}$ as a decimal floating-point number.
\item Comments: This does \emph{not} just call {\tt Log(Log(Float(OrderE8(q))))}, because applying {\tt Float} to too large integers does not result in an appropriate floating-point representation of those integers, but in the output {\tt inf}, which is treated like infinity.
\end{itemize}

\subsubsection{{\tt OrderOutE8}}\label{subsubsec3P17P3}

\begin{itemize}
\item Syntax: {\tt OrderOutE8(q)}, $q$ a prime power.
\item Output: the outer automorphism group order $|\Out(E_8(q))|$.
\end{itemize}

\subsubsection{{\tt NrConjugacyClassesE8}}\label{subsec3P17P4}

\begin{itemize}
\item Syntax: {\tt NrConjugacyClassesE8(q)}, $q$ a prime power.
\item Output: the conjugacy class number $\k(E_8(q))$.
\item Comments: See L{\"u}beck's database \cite{Lue}.
\end{itemize}

\subsubsection{{\tt NrAutOrbitsBoundE8}}\label{subsubsec3P17P5}

\begin{itemize}
\item Syntax: {\tt NrAutOrbitsBoundE8(q)}, $q$ a prime power.
\item Output: the lower bound $\underline{\omega}(E_8(q)):=\lceil\frac{\underline{\k}(E_8(q))}{|\Out(E_8(q))|}\rceil$ on $\omega(E_8(q))$.
\end{itemize}

\subsubsection{{\tt epsilonOmegaBoundE8}}\label{subsubsec3P17P6}

\begin{itemize}
\item Syntax: {\tt epsilonOmegaBoundE8(q)}, $q$ a prime power.
\item Output: the lower bound $\underline{\epsilon_{\omega}}(E_8(q)):=\frac{\log\log{\underline{\omega}(E_8(q))}}{\log\log{|E_8(q)|}}$ on $\epsilon_{\omega}(E_8(q))$.
\end{itemize}

\subsubsection{{\tt CoxeterNoE8}}\label{subsubsec3P17P7}

\begin{itemize}
\item Syntax: {\tt CoxeterNoE8()}.
\item Output: the Coxeter number $h(E_8)=30$.
\end{itemize}

\subsubsection{{\tt SemisimpleElementOrdersE8}}\label{subsubsec3P17P8}

\begin{itemize}
\item Syntax: {\tt SemisimpleElementOrdersE8(q)}, $q=p^f$ a prime power.
\item Output: the set of semisimple (viz., not divisible by $p$) element orders in $E_8(q)$.
\item Comments: Based on \cite{DF91a}.
\end{itemize}

\subsubsection{{\tt NrSemisimpleElementOrdersE8}}\label{subsubsec3P17P9}

\begin{itemize}
\item Syntax: {\tt NrSemisimpleElementOrdersE8(q)}, $q=p^f$ a prime power.
\item Output: the number $\omicron_{\ssrm}(E_8(q))$ of semisimple (viz., not divisible by $p$) element orders in $E_8(q)$.
\item Comments: Simply calls {\tt Length(SemisimpleElementOrdersE8(q))}.
\end{itemize}

\subsubsection{{\tt NrElementOrdersE8Bound}}\label{subsubsec3P17P10}

\begin{itemize}
\item Syntax: {\tt NrElementOrdersE8Bound(q)}, $q=p^f$ a prime power.
\item Output: the upper bound $\overline{\omicron}(E_8(q)):=\omicron_{\ssrm}(E_8(q))\cdot(1+\lceil\log_p(30)\rceil)$ on $\omicron(E_8(q))$, the number of element orders in $E_8(q)$.
\end{itemize}

\subsubsection{{\tt epsilonQBoundE8}}\label{subsubsec3P17P11}

\begin{itemize}
\item Syntax: {\tt epsilonQBoundE8(q)}, $q$ a prime power.
\item Output: the lower bound $\underline{\epsilon_{\q}}(E_8(q)):=\frac{\log\log{(\underline{\omega}(E_8(q))/\overline{\omicron}(E_8(q))+3)}}{\log\log{|E_8(q)|}}$ on $\epsilon_{\q}(E_8(q))$.
\end{itemize}

\subsection{High-end functions for finite simple groups of Lie type}\label{subsec3P18}

\subsubsection{{\tt NrAutOrbits}}\label{subsubsec3P18P1}

\begin{itemize}
\item Syntax: {\tt NrAutOrbits(G)}, $G$ a finite group.
\item Output: the number $\omega(G)$ of $\Aut(G)$-orbits on $G$.
\item Comments: Proceeds by computing the conjugacy classes of $G$ as well as a transversal for $\Inn(G)$ in $\Aut(G)$, and then determining the number of orbits of the action of $\Out(G)$ on the conjugacy classes of $G$.
\end{itemize}

\subsubsection{{\tt ElementOrders}}\label{subsubsec3P18P2}

\begin{itemize}
\item Syntax: {\tt ElementOrders(G)}, $G$ a finite group.
\item Output: the set $\Ord(G)$ of element orders in $G$.
\item Comments: Proceeds by computing first the conjugacy classes of $G$, and then the order in $G$ of a representative of each conjugacy class.
\end{itemize}

\subsubsection{{\tt NrElementOrders}}\label{subsubsec3P18P3}

\begin{itemize}
\item Syntax: {\tt NrElementOrders(G)}, $G$ a finite group.
\item Output: the number $\omicron(G)$ of element orders in $G$.
\item Comments: Simply calls {\tt Length(ElementOrders(G))}.
\end{itemize}

\subsubsection{{\tt epsilonOmegaBoundGeneral1}}\label{subsubsec3P18P4}

\begin{itemize}
\item Syntax: {\tt epsilonOmegaBoundGeneral1(s,Q)}, $s$ one of the strings {\tt "A1"}, {\tt "A2"}, {\tt "2A2"}, {\tt "B2"}, {\tt "C2"}, {\tt "G2"}, {\tt "2B2"} or {\tt "2G2"}, $Q$ a prime power such that if ${^t}X_d$ is the Lie symbol encoded by $s$, then the notation ${^t}X_d(Q)$ is defined.
\item Output: the lower bound on $\epsilon_{\omega}({^t}X_d(Q))$ discussed in \cite[Case (1) of the proof of Theorem 1.1.3(2) in Subsection 3.3]{BGP19a}, based on the data displayed in \cite[Table 3]{BGP19a}.
\end{itemize}

\subsubsection{{\tt epsilonOmegaBoundGeneral2}}\label{subsubsec3P18P5}

\begin{itemize}
\item Syntax: {\tt epsilonOmegaBoundGeneral2(d)}, $d\in\IN^+$, $d\geqslant3$.
\item Output: The GAP evaluation of the expression
\[
\frac{\log(d-\frac{2\log(d+1)}{\log{2}}-\frac{\log{6}}{\log{2}})+\log\log{2}}{\log(4d^2)+\log\log{2}},
\]
which is a uniform lower bound on $\epsilon_{\omega}(S)$ where $S$ is a finite simple Lie type group of the form ${^t}X_d(2^t)$, see \cite[Case (2a) of the proof of Theorem 1.1.3(2) in Subsection 3.3]{BGP19a}.
\item Comments: If $d<10$, the argument in the first logarithm in the numerator is negative, whence the output then is {\tt nan}. For $d\geqslant10$, the output is a decimal floating-point number.
\end{itemize}

\subsubsection{{\tt epsilonOmegaBoundGeneral3}}\label{subsubsec3P18P6}

\begin{itemize}
\item Syntax: {\tt epsilonOmegaBoundGeneral3(d,q)}, $d\in\IN^+$, $d\geqslant3$, $q>2$ either a prime power or a positive algebraic number of the form $\sqrt{2^{2k+1}}$ for some $k\in\IN^+$.
\item Output: the GAP evaluation of the expression
\[
\frac{\log(d-\frac{2\log(d+1)}{\log{q}}-\frac{\log{6}}{\log{q}}-\frac{1}{\e\log{2}})}{\log(4d^2)},
\]
which is a uniform lower bound on $\epsilon_{\omega}(S)$ where $S$ is a finite simple group of Lie type of the form ${^t}X_d(q^t)$, see \cite[Case (2b) of the proof of Theorem 1.1.3(2) in Subsection 3.3]{BGP19a}.
\end{itemize}

\subsubsection{{\tt q0ListOmega}}\label{subsubsec3P18P7}

This is a variable, defined as a certain list of length $16$. For each $d\in\{3,\ldots,18\}$, the list contains a pair of the form {\tt [d,q0]} such that $q_0$ is a prime power and such that for all finite simple Lie type groups $S={^t}X_d(q^t)$ with $q\geqslant q_0$, $\epsilon_{\omega}(S)>\epsilon_{\omega}(\Alt(5))$. See also \cite[Table 4]{BGP19a}.

\subsubsection{{\tt exceptionsListOmega}}\label{subsubsec3P18P8}

\begin{itemize}
\item Syntax: {\tt exceptionsListOmega()}.
\item Output: a finite list, consisting of triples of the form {\tt [s1,d,Q]} and pairs of the form {\tt [s2,Q]} where
\begin{itemize}
\item $s_1$ is one of the strings {\tt "A"}, {\tt "2A"}, {\tt "B"}, {\tt "C"}, {\tt "D"} or {\tt "2D"},
\item $s_2$ is one of the strings {\tt "3D4"}, {\tt "F4"}, {\tt "2F4"}, {\tt "E6"}, {\tt "2E6"}, {\tt "E7"} or {\tt "E8"},
\item $d\in\IN^+$, $d\geqslant3$,
\item $Q$ is a prime power,
\item if ${^t}X$ resp.~${^t}X_d$ is the Lie symbol naturally encoded by $s_1$ resp.~$s_2$, then the notation ${^t}X_d(Q)$ is defined.
\end{itemize}
The corresponding list of finite simple groups of Lie type encoded by the output list has the property that if any nonabelian finite simple group $S$ with $S\not\cong\Alt(5)$ and $\epsilon_{\omega}(S)\leq\epsilon_{\omega}(\Alt(5))$ exists, then $S$ is in that list. See also \cite[end of the proof of Theorem 1.1.3(2) in Subsection 3.3]{BGP19a}.
\end{itemize}

\subsubsection{{\tt nrDistinctPartsPartitions}}\label{subsubsec3P18P9}

This is a variable, defined as a certain list of length $2012$. For each $i\in\{1,\ldots,2012\}$, the $i$-th entry of the list is $s(i-1)$, the number of ordered integer partitions of $i-1$ with pairwise distinct parts. This list of values was copied from the one linked in \cite{OEIS}.

\subsubsection{{\tt g2}}\label{subsubsec3P18P10}

\begin{itemize}
\item Syntax: {\tt g2(d)}, $d\in\IN^+$.
\item Output: the function value $g_2(d):=\sum_{d_++d_-=d}{\sum_{i_+=0}^{d_+}{\sum_{i_-=0}^{d_-}{s(i_+)s(i_-)}}}$ where $s(n)$ denotes the number of ordered integer partitions of $n$ with pairwise distinct parts.
\item Comments: This function is defined and used in \cite[Case (2) of the proof of Theorem 1.1.3(5) in Subsection 3.3]{BGP19a}.
\end{itemize}

\subsubsection{{\tt epsilonQBoundExceptional}}\label{subsubsec3P18P11}

\begin{itemize}
\item Syntax: {\tt epsilonQBoundExceptional(s,Q)}, where
\begin{itemize}
\item $s$ is one of the strings {\tt "2B2"}, {\tt "G2"}, {\tt "2G2"}, {\tt "3D4"}, {\tt "F4"}, {\tt "2F4"}, {\tt "E6"}, {\tt "2E6"}, {\tt "E7"} or {\tt "E8"},
\item $Q$ is a prime power,
\item if ${^t}X_d$ is the Lie symbol naturally encoded by $s$, then the notation ${^t}X_d(Q)$ is defined.
\end{itemize}
\item Output: the GAP evaluation of a certain expression that is a lower bound on $\epsilon_{\q}({^t}X_d(Q))$, defined in \cite[beginning of the proof of Theorem 1.1.3(5) in Subsection 3.3]{BGP19a} and based on the data displayed in \cite[Table 5]{BGP19a}.
\item Comments: As for {\tt epsilonOmegaBoundGeneral2} (Subsubsection \ref{subsubsec3P18P5}), the output may be {\tt nan} for a few (finitely many) pairs $(s,Q)$.
\end{itemize}

\subsubsection{{\tt q0ListQExceptional}}\label{subsubsec3P18P12}

This is a variable, defined as a certain list of length $10$. For each string $s$ of one of the forms {\tt "2B2"}, {\tt "G2"}, {\tt "2G2"}, {\tt "3D4"}, {\tt "F4"}, {\tt "2F4"}, {\tt "E6"}, {\tt "2E6"}, {\tt "E7"} or {\tt "E8"}, the list contains an ordered pair of the form {\tt [s,q0]} where $q_0$ is a prime power such that if ${^t}X_d$ is the (exceptional) Lie symbol naturally encoded by $s$ and $Q$ is a prime power with $Q\geqslant q_0$ such that the notation ${^t}X_d(Q)$ is defined, then $\epsilon_{\q}({^t}X_d(Q))>\epsilon_{\q}(\M)$. See also \cite[Table 5]{BGP19a}.

\subsubsection{{\tt exceptionsListQExceptional()}}\label{subsubsec3P18P13}

\begin{itemize}
\item Syntax: {\tt exceptionsListQExceptional()}
\item Output: a finite list of pairs {\tt [s,Q]} where $s$ is a string of one of the forms {\tt "2B2"}, {\tt "G2"}, {\tt "2G2"}, {\tt "3D4"}, {\tt "F4"}, {\tt "2F4"}, {\tt "E6"}, {\tt "2E6"}, {\tt "E7"} or {\tt "E8"}, and $Q$ is a prime power such that if ${^t}X_d$ is the (exceptional) Lie symbol naturally encoded by $s$, then the notation ${^t}X_d(Q)$ is defined. Moreover, the (finite) list of exceptional finite simple groups of Lie type encoded by the output list has the property that if any exceptional finite simple group of Lie type $S$ with $\epsilon_{\q}(S)\leqslant\epsilon_{\q}(\M)$ exists, then $S$ is in that list.
\end{itemize}

\subsubsection{{\tt epsilonQBoundClassical1}}\label{subsubsec3P18P14}

\begin{itemize}
\item Syntax: {\tt epsilonQBoundClassical1(d,t)}, $d,t\in\IN^+$ ($t$ is interpreted as a \enquote{type} specification).
\item Output: If $t\notin\{1,2,3,4\}$, then it returns the message \enquote{This type is not available. Please set the type to $1$, $2$, $3$ or $4$.} as a string.

If $t=1$, it returns the GAP evaluation of the expression
\[
\frac{\log{((1-\frac{\log{3}}{\log{4}})d-\frac{\frac{2\pi}{\sqrt{3}}\sqrt{d}+3\log{(d+1)}+\log{(2+\frac{\log{2d}}{\log{2}})}+\log{4}}{\log{2}})}+\log\log{2}}{\log{(4d^2)}+\log\log{2}},
\]
which is a uniform lower bound on $\epsilon_{\q}({^t}X_d(2^t))$ for ${^t}X\in\{A,{^2}A,B,C,D,{^2}D\}$; see also \cite[Case (1) of the proof of Theorem 1.1.3(5) in Subsection 3.3]{BGP19a}.

If $t=2$, it returns the GAP evaluation of the expression
\[
\frac{\log{((1-\frac{\log{4}}{\log{9}})d-\frac{\frac{2\pi}{\sqrt{3}}\sqrt{d}+3\log{(d+1)}+\log{(2+\frac{\log{(2d)}}{\log{2}})}+\log{4}}{\log{3}}-\frac{1}{\e\log{2}})}}{\log{(4d^2)}},
\]
which is a uniform lower bound on $\epsilon_{\q}({^t}X_d(q^t))$ for ${^t}X\in\{A,{^2}A,B,C,D,{^2}D\}$ and $q>2$; see also \cite[Case (1) of the proof of Theorem 1.1.3(5) in Subsection 3.3]{BGP19a}.

If $t=3$, it returns the GAP evaluation of the expression
\[
\frac{\log{((1-0.311\cdot\frac{\log{3}}{\log{2}})d-\frac{\log{g_2(d)}+2\log{3}+\log{(2+\frac{\log{(2d)}}{\log{2}})}+\log{2}}{\log{2}})}+\log\log{2}}{\log{(4d^2)}+\log\log{2}},
\]
which is a uniform lower bound on $\epsilon_{\q}({^t}X_d(2^t))$ for ${^t}X\in\{A,{^2}A,B,C,D,{^2}D\}$ if $d\geqslant91$; see also \cite[Case (2) of the proof of Theorem 1.1.3(5) in Subsection 3.3]{BGP19a}.

If $t=4$, it returns the GAP evaluation of the expression
\[
\epsilon_{\q}(S)\geqslant\frac{\log{((1-\frac{\log{4}}{\log{27}})d-\frac{\log{g_2(d)}+2\log{(d+1)}+\log{(2+\frac{\log{(2d)}}{\log{2}})}+\log{2}}{\log{3}}-\frac{1}{\e\log{2}})}}{\log{(4d^2)}},
\]
which is a uniform lower bound on $\epsilon_{\q}({^t}X_d(q^t))$ for ${^t}X\in\{A,{^2}A,B,C,D,{^2}D\}$ if $q>2$ and $d\geqslant55$; see also \cite[Case (2) of the proof of Theorem 1.1.3(5) in Subsection 3.3]{BGP19a}.
\end{itemize}

\subsubsection{{\tt epsilonQBoundClassical2}}\label{subsubsec3P18P15}

\begin{itemize}
\item Syntax: {\tt epsilonQBoundClassical2(d,q)}, $d\in\IN^+$, $q$ a prime power.
\item Output: If $d=1$, it returns the GAP evaluation of the expression
\[
\frac{\log\log{(\frac{q+1}{8\log_2(q)\cdot(\sqrt{(q+1)/2}+\sqrt{(q-1)/2})})}}{\log\log{(q(q^2-1))}},
\]
which is a lower bound on $\epsilon_{\q}(A_1(q))$; see also \cite[Case (4) of the proof of Theorem 1.1.3(5) in Subsection 3.3]{BGP19a}. If $d>1$, it returns the GAP evaluation of the expression
\[
\frac{\log\log{\frac{q^d}{2c(d)\cdot\log_2(q)\cdot\min(d+1,q+1)^2\cdot g_2(d)(1+\lceil\log_2(2d)\rceil)\cdot(q+1)^{d/2}}}}{\log\log{q^{4d^2}}}
\]
where
\[
c(d)=\begin{cases}2, & \text{if }d\not=4, \\ 6, & \text{if }d=4.\end{cases}
\]
Unless the output is {\tt nan}, it is a uniform lower bound on $\epsilon_{\q}({^t}X_d(q^t))$ for ${^t}X\in\{A,{^2}A,B,C,D,{^2}D\}$; see also \cite[Case (4) of the proof of Theorem 1.1.3(5)]{BGP19a}.
\end{itemize}

\subsubsection{{\tt q0ListQClassical}}\label{subsubsec3P18P16}

This is a variable, defined as a certain list of length $53$. For each $d\in\{1,\ldots,53\}$, it contains a pair of the form {\tt [d,q0]} where $q_0$ is a prime power such that for all prime powers $q\geqslant q_0$ and all (classical) Lie symbols ${^t}X\in\{A,{^2}A,B,C,D,{^2}D\}$, $\epsilon_{\q}({^t}X_d(q^t))>\epsilon_{\q}(\M)$. See also \cite[Table 8]{BGP19a}.

\subsubsection{{\tt exceptionsListQClassical}}\label{subsubsec3P18P17}

\begin{itemize}
\item Syntax: {\tt exceptionsListQClassical()}.
\item Output: a finite list of triples {\tt [s,d,Q]} where
\begin{itemize}
\item $s$ is one of the strings {\tt "A"}, {\tt "2A"}, {\tt "B"}, {\tt "C"}, {\tt "D"} or {\tt "2D"},
\item $d\in\IN^+$ and
\item $Q$ is a prime power such that if ${^t}X$ is the (classical) Lie symbol naturally encoded by $s$, then the notation ${^t}X_d(Q)$ is defined.
\end{itemize}
Such a triple {\tt [s,d,Q]} is then to be viewed as a code for the classical finite simple Lie type group ${^t}X_d(Q)$, and the output list encodes a list of such groups. This latter list has the property that if any nonabelian finite simple group $S\not\cong\M$ with $\epsilon_{\q}(S)\leqslant\epsilon_{\q}(\M)$ exists, then $S$ occurs in the list. See also \cite[Case (5) of the proof of Theorem 1.1.3(5)]{BGP19a}.
\end{itemize}

\end{document}